# Location and Capacity Planning of Facilities with General Service-Time Distributions Using Conic Optimization

Amir Ahmadi-Javid[†], Oded Berman[*], Pooya Hoseinpour[‡]

[†]*Department of Industrial Engineering & Management Systems, Amirkabir University of Technology, Tehran, Iran*
[*]*Rotman School of Management, University of Toronto, Toronto, Canada*
[‡]*Desautels Faculty of Management, McGill University, Montreal, Quebec, Canada*

**Abstract.** This paper studies a stochastic congested location problem in the network of a service system that consists of facilities to be established in a finite number of candidate locations. Population zones allocated to each open service facility together creates a stream of demand that follows a Poisson process and may cause congestion at the facility. The service time at each facility is stochastic and depends on the service capacity and follows a general distribution that can differ for each facility. The service capacity is selected from a given (bounded or unbounded) interval. The objective of our problem is to optimize a balanced performance measure that compromises between facility-induced and customer-related costs. Service times are represented by a flexible location-scale stochastic model. The problem is formulated using quadratic conic optimization. Valid inequalities and a cut-generation procedure are developed to increase computational efficiency. A comprehensive numerical study is carried out to show the efficiency and effectiveness of the solution procedure. Moreover, our numerical experiments using real data of a preventive healthcare system in Toronto show that the optimal service network configuration is highly sensitive to the service-time distribution. Our method for convexifying the waiting-time formulas of M/G/1 queues is general and extends the existing convexity results in queueing theory such that they can be used in optimization problems where the service rates are continuous.

**Keywords:** Service system design; Stochastic location models; Congested networks; M\M\1 and M\G\1 queue systems; Location problems with congestion; Integer nonlinear programming; Second-order cone programming (Conic quadratic programming or Mixed 0-1 conic optimization); Mixed-integer convex optimization.

# 1 Introduction

This paper focuses on the design of congested service facilities with immobile servers. In such systems congestion arises because population zones (customers) generate streams of stochastic demand for service, which evolve according to Poisson processes, the resources of the facilities are limited, and the service times provided by the facilities are stochastic. The design problem determines the optimal number, locations, and service capacities of the established facilities. Applications of the design problem include private service facilities such as retail stores or maintenance shops, and public service facilities, such as government offices and hospitals.

In a survey on Stochastic Location models with Congestion and Immobile Servers (SLCIS models), Berman and Krass (2015) classified SLCIS models into four classes according to the type of demand for service: inelastic or elastic, and according to the way demand is assigned to facilities which is determined by: a central authority (a decision maker) or by the customers (user choice). In the models presented in our paper, the demand is assumed to be inelastic and allocation of demand can be assumed to be determined either by a decision maker or by the customers. The paper also belongs to the family of balanced-objective models. The models that belong to this family are designed to balance the travel and congestion costs that are incurred by the customers, and the facility-related costs that include the fixed costs of opening facilities and the costs to serve customers.

Balance-objective SLCIS models are studied under different assumptions in many papers, such as Marianov and Rios (2000), Wang et al. (2004), Elhedhli (2006), Berman and Drezner (2007), Aboolian et al. (2008), Castillo et al. (2009), Abouee-Mehrizi et al. (2011), Vidyarthi and Jayaswal (2014), Ahmadi-Javid and Hoseinpour (2017). The main difficulty of balance-objective SLCIS models stems from the nonlinearity of the congestion term.

Except for a few papers, e.g., Berman and Drezner (2007), and Aboolian et al. (2008), which considered M/M/k queue systems in their SLCIS models, each facility is modeled via a one-server queue, and its service capacity is treated as a decision variable. The advantage of considering one-server facilities is that it is easier to handle mathematically and is more practical when the facility uses a variety of processing resources (Berman & Krass, 2015).

There are two approaches for determining the service capacities for one-server facilities. One approximate approach, which we call the *discretization* approach, is to assume that there is a finite set of potential service capacities and one of them should be selected. The other approach, which we call the *continuous* approach, is to assume that an interval of service capacities is given and a point on this interval should be selected. The continuous approach is the more desirable one since it can cover the real-world cases where the intervals can be long or unbounded.

The continuous approach has been broadly studied in the queuing literature; see, for example, Weber (1983), Grassmann (1983), Tu and Kumin (1983), Harel and Zipkin (1987), Shaked and Shanthikumar (1988), Buzacott and Shanthikumar (1992), Fridgeirsdottir and Chiu (2005). However, the majority of the work on SLCLS models is based on the discretizing approach because of its simplicity. Recently, a couple of papers used the continuous approach only for M/M/1 queues. Castillo et al. (2009) and Hoseinpour and Ahmadi-Javid (2016) applied this approach to develop balance-objective SLCLS models, and presented Lagrangian heuristics to solve their models. Aboolian et al. (2015) used the continuous approach in a very different setting where demands are elastic. They solved their model by linearizing the nonlinear terms in the equilibrium constraints.

In most of the SLCLS models, it is assumed that the distribution of the service time is exponential due to its simplicity. This assumption may not be realistic in practice (Boffey et al., 2006). For the M/M/1 system, since the variance of the service-time distribution is equal to the square of the mean service time (the reciprocal of the service capacity), knowledge of the mean value is sufficient. However, for the M/G/1 system, the variance of the service time appears in the expected waiting time in the system, which is the key performance measure used in the SLCLS context. As discussed in Berman and Krass (2015), quite often it is required to determine not only how much capacity to add, but also what kind of capacity to install. For example, the variance of the service time may be small when a higher automation is added and may be large when a manual capacity is added, which can result in the same mean service time but different distributions, and different variances in particular.

Recently, Vidyarthi and Jayaswal (2014) and Ahmadi-Javid and Hoseinpour (2017) applied the discretization approach to study balance-objective SLCLS models with M/G/1 queues. The input of these models for each facility is a finite set of pairs of service capacities and service-time variances, from which one pair should be selected. In these models, the variance is treated as an intrinsic parameter, which is not affected by the service capacity. However, under the continuous approach it is required to allow the variance of the service-time distribution to depend on the service capacity, which makes the resulting model very complicated.

In our paper, for the first time a balance-objective SLCLS model with M/G/1 facilities under the continuous approach is considered, where each service-time variance is a function of the service capacity. To determine how the expected waiting time in an M/G/1 queue changes with the service capacity, we represent the service time of each facility as a general location-scale stochastic model, which covers many well-known univariate distributions. This model extends the simple scale model used in the queuing literature to develop convexity results for G/G/1 queues by Weber (1983), Stidham (1992), and Fridgeirsdottir and Chiu (2005).

To solve the proposed problem, we first formulate it as a mixed-integer non-convex minimization model. Then, we reformulate the model as a Mixed-Integer Second-Order Cone Program (MISOCP), which makes it possible to solve the problem efficiently. Valid inequalities and a cut-generation procedure are also used to enhance the solution procedure in solving large-size instances. Recently, Ahmadi-Javid and Hoseinpour (2017) proposed MISOCPs for an SLCLS problem under the discrete approach where a finite set of candidate pairs of service rate and service-time variance is given as input for each service-time distribution. They showed that the convexified models significantly outperform the best-known solution method presented for that problem. However, their reformulation cannot be adapted for our new SLCLS problem, in which service rates are continuous and the variances are nonlinear functions of these variables (see Remark 5 in Section 7). Moreover, the results obtained by Weber (1983), which show that the queuing formulas are convex in continuous service rates, is not helpful here because i) the underling stochastic representation model (see (23)) is very simple and does not cover well-known distributions, and ii) the convexity of the formulas does not hold jointly in both arrival and service rates even for the M/M/1 queues (see Remark 8 in Section 7), and hence the resulting formulation is not convex.

The main contributions of this paper are as follows:

1- For the first time, this paper studies an SLCLS problem with M/G/1 queues where service rates are continuous decision variables. The resulting optimization model is a mixed-integer nonlinear program that can be solved optimally by our solution method developed based on mixed-integer convex optimization.

2- To express service-time distributions in terms of service rates, a flexible stochastic representation model is used, which covers many known distributions and considerably generalizes the simple scale model in queuing theory. Our convexification method for M/G/1 queues extends the convexity results developed for G/G/1 queues with the scale representation model in a way that the resulting large-scale optimization models can be solved optimally by existing optimization solvers.

3- MISOPCs are promising mixed-integer convex optimization models that can be solved in large scales by the state-of-the-art optimization solvers such as CPLEX. The condition under which the model of our SLCLS problem can be cast as an MISOCP is investigated (Theorem 1). This condition is proven to hold for the proposed stochastic representation model (Proposition 1). Moreover, more efficient MISOCPs are obtained for special cases (Section 4). These reformulations enable us to solve the problem exactly. In fact, our reformulations extend the domain of application of MISOCPs. Our method is general and can be applied to similar problems (Section 7).

4- A class of valid inequities is proposed to enhance the MISOCPs (Section 5.1). A branch-and-cut algorithm is developed to solve the problem, in which effective *lift-and-project* cuts are iteratively added (Section 5.2). The numerical results indicate the effectiveness and efficiency of our solution method to solve large-size problem instances (Section 6.1). This is an important achievement because a much simpler version of our problem (the problem with M/M/1 queues and without any upper bounds on service rates) has been solved only heuristically.

5- Our study contributes to the service operations management by highlighting that the service process and network configuration in a service system can be dependent. As an example, our numerical experiments based on real data of the city of Toronto indicate that the network configuration of an immobile service system can be significantly different when changing the service time distribution (Section 6.2). More specifically we obtain very different configurations when we change the service time distribution from exponential to Gamma.

The remainder of the paper is organized as follows. Section 2 states and mathematically models the problem. Section 3 provides our main results, which show that the resulting model can be reformulated as an MISOCP. Section 4 investigates special cases for which simpler MISOCPs can be obtained. Section 5 provides details on the branch-and-cut algorithm used to solve the proposed MISOCPs. Section 6 reports the computational results and the application of the model to an example based on real data. Section 7 extends the reformulation method and provides general convex optimization results for waiting-time metrics of M/G/1 queues. Section 8 concludes the paper and includes suggestions for future studies. Appendix A provides some preliminaries on MISOCPs, and Appendix B provides the proof of Proposition 4. Appendix C presents an alternative MISOCP for the problem.

## 2   Problem statement and mathematical modeling

We call the problem considered in this paper the Service-System Design Problem (SSDP). In this problem, given a network of potential locations for installing immobile Service Facilities (SFs), a set of SFs are selected to establish for providing quality service to the customers. Customers are geographically classified to demand zones. The SSDP simultaneously determines locations of open SFs, their service capacities together with the allocation of demand zones to the open SFs in order to balance the system's service cost and total customer cost. The service cost includes the costs of establishing SFs and of providing service. The total customer cost is the sum of customers' travel and waiting costs, where the smaller the sum, the better the service quality. The stochastic nature of customers' demand along with uncertain service time may cause queues of customers waiting for service in the SFs, which are modelled as M/G/1 queuing systems.

Let $I$ be the set of candidate SFs and $J$ be the set of demand zones, then define the following decision variables:

- $x_i$    A binary variable that equals 1 if SF $i \in I$ is established, and 0 otherwise.
- $y_{ij}$    A binary variable that equals 1 if demand zone $j \in J$ is assigned to SF $i \in I$ to be served, and 0 otherwise.
- $\mu_i$    A non-negative real-valued variable that represents the service capacity (rate) at SF $i \in I$.

The service time at SF $i$ is represented by a nonnegative random variable $S_i$ with finite mean $E(S_i) = 1/\mu_i$ and variance $var(S_i)$, where the service rate $\mu_i$ is a decision variable. Based on the distribution of $S_i$, $var(S_i)$ is assumed to depend on the service rate $\mu_i$. Therefore, it is simply represented by a function of $\mu_i$ as $var(S_i) = v_i(\mu_i)$ where $v_i: \mathbb{R}_{\geq 0} \to \mathbb{R}_{\geq 0}$ is called the variance function for service time $S_i$. For example, when $S_i$ follows the exponential distribution $\text{Exp}(1/\mu_i)$, then we have $v_i(\mu_i) = \mu_i^{-2}$, while for $S_i \sim \text{Uniform}(1/\mu_i - \theta_i, 1/\mu_i + \theta_i)$ with $\theta_i$ being a positive constant, the variance function becomes a constant, i.e., $v_i(\mu_i) = \theta_i^2/3$.

We note that the distribution of $S_i$ can take only nonnegative values. To ensure it, the set of feasible values for decision variable $\mu_i$ must be adjusted appropriately. To keep $S_i$ a nonnegative random variable, three cases for the distribution of $S_i$ are possible from a modeling point of view as follows:

**Case 1**. The essential infimum of the distribution of $S_i$ is finite and independent of $\mu_i$. In this case, if the essential infimum is greater than or equal to zero, for any positive value of $\mu_i$, the random variable $S_i$ is nonnegative. The Exponential and Gamma distributions are instances of such distributions, whose essential infimums are always zero.

**Case 2**. The essential infimum of the distribution of $S_i$ is finite, and depends on $\mu_i$. In this case, only those positive values of $\mu_i$ are acceptable for which the essential infimum is greater than or equal to zero. For example, if $S_i$ is uniformly distributed over the interval $(1/\mu_i - \theta_i, 1/\mu_i + \theta_i)$ with $\theta_i > 0$, the service rate $\mu_i$ must satisfy $0 < \mu_i \leq 1/\theta_i$ to guarantee that $S_i$ is a nonnegative random variable.

**Case 3**. The essential infimum of the distribution of $S_i$ is negative infinity. This may be considered when a distribution with infinite essential infimum is not appropriate to model the distribution of the service time $S_i$, but by a little adjustment, if possible, such distribution can also be used for this purpose. For example, it is possible to only consider those values of the service rate $\mu_i$ for which the probability of $S_i$ being negative is negligible as follows:

$$\Pr\{S_i < 0\} \leq \alpha \text{ or } \Pr\{S_i \geq 0\} \geq 1 - \alpha,$$

where the probability $\alpha \in (0,1)$ should be set to a very small number, such as $10^{-2}$. For example, when the service time follows a normal distribution, $S_i \sim \text{Normal}(1/\mu_i, \sigma_i^2)$, from the above inequality one obtains $0 < \mu_i \leq 1/z_\alpha \sigma_i$, with $z_\alpha$ the $\alpha$ upper percentile of the standard normal distribution.

The set of possible values of service rate $\mu_i$ is represented by $m_i \leq \mu_i \leq M_i, \mu_i > 0$ where the bounding parameters $0 < m_i \leq M_i$ must be set such that the random variable $S_i$ is always (or with high probability) nonnegative. Besides this modeling advantage, these bounds can be used to restrict the feasible values for the service rate $\mu_i$ to those that can be implemented in practice based on the decision maker's preference.

The demand of each customer $j$ evolves according to a Poisson process with intensity rate $\lambda_j > 0$. Hence, using the superposition property of Poisson processes, the aggregated demand from SF $i$ also follows a Poisson process with rate $\sum_{j \in J} \lambda_j y_{ij}$. Then, for an M/G/1 queue, the expected sum of the customers' waiting time at SF $i$, denoted by a bivariate function $WT_i(\cdot,\cdot)$, is given by

$$WT_i\left(\sum_{j \in J} \lambda_j y_{ij}, \mu_i\right) = \left(\sum_{j \in J} \lambda_j y_{ij}\right)\left(\frac{\left(\sum_{j \in J} \lambda_j y_{ij}\right)(1 + v_i(\mu_i)\mu_i^2)}{2\mu_i(\mu_i - \sum_{j \in J} \lambda_j y_{ij})} + \frac{1}{\mu_i}\right)$$

whenever the steady-state condition $\sum_{j \in J} \lambda_j y_{ij} < \mu_i$ holds (Gross, 2008). Therefore, using the following cost parameters:

- $ec_i$     Fixed cost per time unit for establishing SF $i \in I$
- $wc_i$     Waiting cost per time unit for a customer at SF $i \in I$,
- $sc_i$     Service provision cost per service at SF $i \in I$
- $tc_{ij}$     Transportation cost per travel between customer $j \in J$ and SF $i \in I$

the proposed SSDP can be formulated as follows:

$$\min \sum_{i \in I} ec_i x_i + \sum_{i \in I} sc_i \mu_i + \sum_{i \in I} wc_i WT_i\left(\sum_{j \in J} \lambda_j y_{ij}, \mu_i\right) + \sum_{i \in I}\sum_{j \in J} tc_{ij} \lambda_j y_{ij} \quad (1)$$

s.t.

$$\sum_{i \in I} y_{ij} = 1 \qquad j \in J \quad (2)$$

$$y_{ij} \leq x_i \qquad i \in I, j \in J \quad (3)$$

$$\sum_{j \in J} \lambda_j y_{ij} \leq \mu_i \qquad i \in I \quad (4)$$

$$m_i x_i \leq \mu_i \leq M_i x_i \qquad i \in I \quad (5)$$

$$x_i \in \{0,1\}, y_{ij} \in \{0,1\} \qquad i \in I, j \in J. \quad (6)$$

Objective function (1) includes the service system's operational costs (the first and second terms) and the customers' costs (the third and fourth terms). The first term is the total cost of establishing SFs. The second term is the expected total cost of providing service at open SFs. The third term is the sum of expected customers' waiting costs, while the fourth term is the sum of expected customers' travel costs. To avoid undetermined cases in the objective function, we accept the conventions $\frac{a}{0} = +\infty$ for $a > 0$ and $\frac{0}{0} = 0$. One may want to consider different weights for these cost terms; however, for notational simplicity, without any loss of generality, here it is assumed all terms be equally weighted to one. Indeed, the impact of weights can directly be incorporated into the cost parameters.

Constraints (2) ensure that each customer is assigned to only one open SF. Constraints (3) prevent a customer from being assigned to a non-established SF. Constraints (4) force the required steady-state conditions. Constraints (5) impose the required lower and upper bounds on the service rates of open SFs. Constraints (6) include binary decision variables.

Here, for simplicity, we assumed that the service-provision cost for each SF linearly depends on the service capacity (which is also commonly used in the literature). However, any other alternative cost structure, such as piecewise-linear functions, that are (exactly or approximately) representable by a mixed-integer linear system could be used. Moreover, other linear constraints can also be added to the model (1)-(6). For example, to ensure that each customer is assigned to the closest open SF, the following constraints can be added:

$$\sum_{k \in I} y_{kj} d_{kj} + (D_i - d_{ij}) x_i \leq D_i \qquad i \in I, j \in J \qquad (7)$$

where $d_{ij}$ denotes the distance between customer $j$ and SF $i$, and where $D_i$ is a big positive constant satisfying $D_i \geq \max_{j \in J}\{d_{ij}\}$. Other alternatives for constraints (7) are presented and compared in Gerrard and Church (1996) and Berman et al. (2006). It is important to note that the formulation without constraints (7) is suitable to central-authority assignment, and with constraints (7) is suitable to user-choice assignments.

The model (1)–(6) is a mixed-integer nonlinear program with non-convex minimization objective and linear constraints. Solving large-size instances of such problems is generally challenging (Köppe, 2012). In the next sections, it is shown that for a broad class of service-time distributions, the model (1)–(6) can be reformulated as MISOCPs, which are efficiently solvable in large scales.

MISOCPs have recently been applied for modeling various problems in portfolio optimization (Benson and Saglam, 2013), options pricing (Pinar, 2013), telecommunication network design (Hijazi et al., 2013), power distribution system planning (Taylor and Hover, 2012), battery swapping stations, (Mak et al., 2013), Euclidean k-center, (Brandenberg and Roth, 2009), supply chain network design (Ahmadi-Javid & Azad, 2010; Atamturk et al., 2012), and berth allocation (Du et al., 2011). For an overview of available applications see Benson and Saglam (2013). More technical details and references on algorithms developed for solving MISOCPs are provided in Section 5.

# 3  MISOCP reformulation

As stated in Section 2, the model (1)–(6) is a complex mixed-integer non-convex program. This section shows how this model can be reformulated as MISOCPs, which are known to be efficiently solvable mixed-integer convex programs. For a brief introduction to MISOCPs see Appendix A. The structure of the objective function (1) depends on the variance functions $v_i(.), i \in I$. Theorem 1 provides conditions for the variance functions under which the non-convex model (1)–(6) can be convexified as an MISOCP (see Appendix A for basic definitions).

**Theorem 1.** Let $\sigma_i, i \in I$ be nonnegative variables. If the constraints

$$v_i(\mu_i) \le \sigma_i^2 \qquad i \in I, \qquad (8)$$

are SOC-representable (i.e., they can be rewritten by a set of SOC constraints), then the model (1)–(6) can be reformulated as the following MISOCP:

$$\min \sum_{i \in I} \left( ec_i x_i + sc_i \mu_i + wc_i(t_i + \rho_i) + \sum_{j \in J} tc_{ij} \lambda_j y_{ij} \right) \qquad (9)$$

s.t.

(2)–(6), (8)

$$\sigma_i - (1 - y_{ij}) Q_i \le u_{ij} \le \sigma_i \qquad i \in I, j \in J \qquad (10)$$

$$0 \le u_{ij} \le y_{ij} Q_i \qquad i \in I, j \in J \qquad (11)$$

$$\left\| \begin{matrix} 2\sqrt{\lambda_1} y_{i1} \\ \vdots \\ 2\sqrt{\lambda_{|J|}} y_{i|J|} \\ \rho_i - \mu_i \end{matrix} \right\| \le \rho_i + \mu_i \qquad i \in I \qquad (12)$$

$$\left\| \begin{matrix} 2\rho_i \\ 2\sum_{j \in J} \lambda_j u_{ij} \\ 2 - 2\rho_i - t_i \end{matrix} \right\| \le 2 - 2\rho_i + t_i \qquad i \in I \qquad (13)$$

$$t_i \ge 0, \rho_i \ge 0, \qquad i \in I \qquad (14)$$

where $Q_i, i \in I$, is a constant upper bound for $\sigma_i$ that must satisfy

$$\sup_{m_i' \le \mu_i \le M_i} \{v_i(\mu_i)\} < Q_i < +\infty \qquad i \in I \qquad (15)$$

with $m_i' = \max \left\{ \min_{j \in J} \{\lambda_j\}, m_i \right\}$.

**Proof.** By introducing variable $\rho_i \ge 0$ for each $i \in I$ such that

$$\sum_{j \in J} \lambda_j y_{ij} = \rho_i \mu_i \qquad i \in I, \qquad (16)$$

the objective (1) can be transformed to the following one:

$$\min \sum_{i \in I} \left( ec_i x_i + sc_i \mu_i + wc_i \left( \frac{\rho_i^2 + \left(\sum_{j \in J} \lambda_j y_{ij}\right)^2 \sigma_i^2}{2(1 - \rho_i)} + \rho_i \right) + \sum_{j \in J} tc_{ij} \lambda_j y_{ij} \right). \tag{17}$$

By defining new variables $t_i \geq 0$ for each $i \in I$ such that

$$\rho_i^2 + \left(\sum_{j \in J} \lambda_j u_{ij}\right)^2 = 2(1 - \rho_i) t_i \qquad i \in I \tag{18}$$

with

$$u_{ij} = \sigma_i y_{ij} \qquad i \in I, j \in J, \tag{19}$$

the model becomes

$$\min \sum_{i \in I} \left( ec_i x_i + sc_i \mu_i + wc_i (t_i + \rho_i) + \sum_{j \in J} tc_{ij} \lambda_j y_{ij} \right) \tag{20}$$

s.t.
(2)–(6), (8), (16), (18), (19).

Objective (17) minimizes a function that includes the term

$$\frac{\rho_i^2 + \left(\sum_{j \in J} \lambda_j y_{ij}\right)^2 \sigma_i^2}{2(1 - \rho_i)} + \rho_i,$$

which is increasing in $\rho_i$. Thus, because $\rho_i$ appears only in (17), and noting that $y_{ij} = y_{ij}^2$ for binary variables $y_{ij} \in \{0,1\}$, one can replace identity (16) by inequality (21) below for each $i \in I$.

$$\sum_{j \in J} \lambda_j y_{ij}^2 \leq \rho_i \mu_i \qquad i \in I. \tag{21}$$

Also, constraint (18) can be replaced by (22) because the variable $t_i$ appears in the minimization objective (20).

$$\rho_i^2 + \left(\sum_{j \in J} \lambda_j u_{ij}\right)^2 \leq 2(1 - \rho_i) t_i \qquad i \in I. \tag{22}$$

For appropriate large positive numbers $Q_i$ satisfying (15), constraints (19) can be linearized as constraints (10) and (11). Constraints (21) and (22) are hyperbolic constraints, which can be represented by SOC constraints (12) and (13), respectively. This completes the proof. ∎

Theorem 1 indicates the important role of constraints (8) in reformulating our problem as an MISOCP. These constraints depend on the distributions of the service times through the variance functions. Therefore, to continue our analysis, it is required to determine how service times depend on service capacities, which reveals the structures of the variance functions.

Unlike the M/M/1 systems, the arrival rate and service capacity are not necessarily sufficient to completely describe performance measures of M/G/1 queues. Therefore, a method is required to express how the service time in an M/G/1 queue depend on its service capacity (Fridgeirsdottir and Chiu, 2005). A common method used in research studies on M/G/1 and G/G/1 queues is to represent each service time as a simple function of the service capacity and a random variable (Fridgeirsdottir & Chiu, 2005; Stidham, 1992; Weber, 1983). Actually, these studies assume that service time $S$ at each SF depends on service capacity $\mu$ by the following scale model:

$$S = \frac{1}{\mu}X, \tag{23}$$

where $X$ is a positive random variable with $E(X) = 1$ and with finite variance. Extending this method, our paper considers a more general, flexible location-scale model for representing service times. For each SF $i \in I$, let each service time $S_i$ be represented as follows:

$$S_i = \frac{1}{\mu_i} + \sum_{l=0}^{L_i} \left(\frac{1}{\mu_i}\right)^l \delta_{i,l} \varepsilon_{i,l} \qquad i \in I, \tag{24}$$

where $\varepsilon_{i,l}, l = 0, \cdots, L_i$, are independent random variables with $E(\varepsilon_{i,l}) = 0$ and $var(\varepsilon_{i,l}) = 1$. The parameters $\delta_{i,l}, l = 0, \cdots, L_i$, are positive constants, which are sort of scale parameters while $1/\mu_i$ is a location parameter, which can also impact the dispersion of the service time $S_i$ whenever $L_i \geq 1$.

By setting $L_i = 1$, $\delta_{i,l} = 1$, $\varepsilon_{i,0} = 0$, and assuming essinf $\varepsilon_{i,1} = -1$, one can retrieve from (24) the traditional scale model given in (23) as follows:

$$S_i = \frac{1}{\mu_i} + \varepsilon_{i,0} + \frac{1}{\mu_i}\varepsilon_{i,1} = \frac{1}{\mu_i}(1 + \varepsilon_{i,1}),$$

where $1 + \varepsilon_{i,1}$ is essentially the same as $X$ in (23).

From (24), the mean and variance of each service time are given by

$$E(S_i) = \frac{1}{\mu_i}, var(S_i) = v_i(\mu_i) = \sum_{l=0}^{L_i} \left(\frac{1}{\mu_i}\right)^{2l} \delta_{i,l}^2 \qquad i \in I. \tag{25}$$

Therefore, under the setting in (25) constraints (8) become as follows:

$$\sum_{l=0}^{L_i} \left(\frac{1}{\mu_i}\right)^{2l} \delta_{i,l}^2 \leq \sigma_i^2 \qquad i \in I. \tag{26}$$

Proposition 1 below shows that these constraints are SOC-representable, which implies that Theorem 1 can be applied whenever service times follow the flexible structure given in (24). The following lemma is required to prove the forthcoming proposition.

**Lemma 1.** For any given integer $p$, constraint $y^p \leq t$ is SOC-representable.

**Proof.** Along the same line used in Ben-Tal and Nemirovski (2001), one can rewrite the constraint in the following form:

$$y \leq (z_1 z_2 \cdots z_{2^l})^{\frac{1}{2^l}}$$

$$z_1, z_2, \cdots, z_r = y, z_{r+1} = t, z_{r+2}, \ldots, z_{2^l} = 1,$$

where $l$ is the smallest integer satisfying $p \leq 2^l$ and $r = 2^l - p$. Now $y \leq (z_1 z_2 \cdots z_{2^l})^{1/2^l}$ can be reformulated as follows:

$$
\begin{aligned}
& z_i = z_{0,i} & & i = 1, \ldots, 2^l \\
& z_{1,i} \leq \sqrt{z_{0,2i-1} z_{0,2i}}, \quad z_{1,i}, z_{0,2i-1}, z_{0,2i} \geq 0 & & i = 1, \ldots, 2^{l-1} \\
& z_{2,i} \leq \sqrt{z_{1,2i-1} z_{1,2i}}, \quad z_{2,i}, z_{1,2i-1}, z_{1,2i} \geq 0 & & i = 1, \ldots, 2^{l-2} \\
& \quad \vdots & & \\
& z_{l,i} \leq \sqrt{z_{l-1,1} z_{l-1,2}}, \quad z_{l,1}, z_{l-1,1}, z_{l-1,2} \geq 0 & & \\
& s \leq z_{l,1}. & &
\end{aligned}
$$

Each hyperbolic constraint $u \leq \sqrt{vw}$ with $u, v, w \geq 0$ is identical to the following SOC constraint:

$$\sqrt{4u^2 + (v-w)^2} \leq v + w.$$

This completes the proof. ∎

**Proposition 1.** If service times can be represented in the form of (24), then constraints (8), which are already refined in (26), are SOC-representable.

**Proof.** Constraints (26) can be rewritten as

$$z_i \geq \frac{1}{\mu_i} \qquad i \in I \tag{27}$$

$$\sigma_i^2 \geq \delta_{i,0}^2 + \sum_{l=1}^{L_i} \delta_{i,l}^2 z_i^{2l} \qquad i \in I \tag{28}$$

where $z_i, i \in I$, are new nonnegative variables. The above constraints can be represented by

$$\sqrt{4 + (z_i - \mu_i)^2} \leq z_i + \mu_i \qquad i \in I \tag{29}$$

$$\sigma_i \geq \sqrt{\delta_{i,0}^2 + \sum_{l=1}^{L_i} \delta_{i,l}^2 \bar{z}_{i,l}^2} \qquad i \in I \tag{30}$$

$$z_i^l \leq \bar{z}_{i,l} \qquad i \in I, l = 2, \ldots, L_i, \tag{31}$$

where $\bar{z}_{i,l}, i \in I, l = 2, \ldots, L_i$, are new nonnegative variables with $\bar{z}_{i,0} = 1$ and $\bar{z}_{i,1} = z_i$. Constraints (29) and (30) are already SOC constraints, and constraints (31) are also SOC-representable by Lemma 1. This ends the proof. ∎

From a queuing-theory perspective, one can deduce that Theorem 1 significantly extends the old convexity results obtained by for the simple scale model (see (23)) for G/G/1 queues by Weber (1983) and Stidham (1992) in a way that the resulting model can be solved by existing optimization solvers.

Note that the same result in Proposition 1 can be obtained for the following more general representation for each $i \in I$:

$$S_i = \frac{1}{\mu_i} + \sum_{l=-L_i^-}^{L_i^+} \left(\frac{1}{\mu_i}\right)^l \delta_{i,l} \varepsilon_{i,l}$$

where are $L_i^-$ and $L_i^+$ are non-negative integers by Lemma 1. However, this model is of less interest as the associated variance can positively depend on service rate.

## 4  Simplified MISOCPs for special cases

In Section 3, the general form of the service times presented by (24) is considered. In this section, we analyze some special classes of service-time distributions for which simpler forms of MISOCP with a smaller number of variables and constraints can be proposed.

For the special case of $L_i = 0$, one retrieves the classic location-scale model

$$S_i = \frac{1}{\mu_i} + \delta_{i,0}\varepsilon_{i,0} \qquad i \in I, \tag{32}$$

with $var(S_i) = \delta_{i,0}^2$ that is independent of the service rate $\mu_i$. This case is referred to as the constant case in the sequel.

The affine case is given by

$$S_i = \frac{1}{\mu_i} + \delta_{i,0}\varepsilon_{i,0} + \frac{1}{\mu_i}\delta_{i,1}\varepsilon_{i,1} \qquad i \in I, \tag{33}$$

which is obtained from (24) by setting $L_i = 1$. For the affine case,

$$var(S_i) = \delta_{i,0}^2 + \delta_{i,1}^2/\mu_i^2.$$

These constant and affine cases cover most practical cases and many important distributions, such as exponential, double exponential (Laplace), gamma, normal, exponential power, stable, and lognormal distributions.

**Proposition 2.** When all service times follow the affine case, given in (33), the model (1)–(6) can be rewritten as the following MISOCP model:

$$\min \sum_{i \in I} \left( ec_i x_i + sc_i \mu_i + wc_i(t_i + \rho_i) + \sum_{j \in J} tc_{ij} \lambda_j y_{ij} \right) \tag{34}$$

s.t.
(2)–(6), (12), (14)

$$\left\| \begin{array}{c} 2(1+\delta_{i,1}^2)^{1/2} \rho_i \\ 2\delta_{i,0} \sum_{j \in J} \lambda_j y_{ij} \\ 2 - 2\rho_i - t_i \end{array} \right\| \leq 2 - 2\rho_i + t_i \qquad i \in I. \tag{35}$$

**Proof.** Substituting $\sigma_i^2 = \delta_{i,0}^2 + \delta_{i,1}^2/\mu_i^2$ into the objective function (1), the model (1)–(6) can be reformulated as

$$\min \sum_{i \in I} \left( ec_i x_i + sc_i \mu_i + wc_i(t_i + \rho_i) + \sum_{j \in J} tc_{ij} \lambda_j y_{ij} \right) \tag{36}$$

s.t.
(2)–(6), (12), (14)

$$(1+\delta_{i,1}^2)\rho_i^2 + \left( \delta_{i,0} \sum_{j \in J} \lambda_j y_{ij} \right)^2 \leq 2(1-\rho_i)t_i \qquad i \in I, \tag{37}$$

where (37) can be represented as SOC (35). This completes the proof. ∎

**Proposition 3.** When for all the service-time distributions it holds that $var(S_i) = 1/\mu_i^2$, $m_i = 0$, and $M_i \to +\infty$ (cf. (25)), which includes the important case of M\M\1, the model (1)–(6) can be rewritten as the following MISOCP model:

$$\min \sum_{i \in I} ec_i x_i + \sum_{i \in I} \sum_{j \in J} sc_i \lambda_j y_{ij} + \sum_{i \in I} r_i + \sum_{i \in I} \sum_{j \in J} tc_i \lambda_j y_{ij} \tag{38}$$

s.t.
(2)–(3)

$$\sqrt{wc_i sc_i \sum_{j \in J} \lambda_j y_{ij}^2} \leq r_i \qquad i \in I \tag{39}$$

$$r_i \geq 0 \qquad i \in I \tag{40}$$

$$x_i \in \{0,1\}, y_{ij} \in \{0,1\} \qquad i \in I, j \in J. \tag{41}$$

**Proof.** For the affine case (33) with $\delta_{i,1}^2 = 1, \delta_{i,0}^2 = 0$, when $m_i = 0$ and $M_i \to +\infty$, the proposed model (1)–(6) can be elaborated as follows:

$$\min \sum_{i \in I} ec_i x_i + \sum_{i \in I} sc_i \mu_i + \sum_{i \in I} \frac{wc_i \sum_{j \in J} \lambda_j y_{ij}}{\mu_i - \sum_{j \in J} \lambda_j y_{ij}} + \sum_{i \in I} \sum_{j \in J} tc_i \lambda_j y_{ij} \tag{42}$$

s.t.
(2)–(4), (6).

Thus, the optimal value of $\mu_i$ can explicitly be found as

$$\mu_i^* = \sqrt{\frac{wc_i \sum_{j \in J} \lambda_j y_{ij}}{sc_i}} + \sum_{j \in J} \lambda_j y_{ij}.$$

Substituting $\mu_i^*$ into the objective function (36) results in

$$\min \sum_{i \in I} ec_i x_i + \sum_{i \in I} \sqrt{sc_i wc_i \sum_{j \in J} \lambda_j y_{ij}} + \sum_{i \in I} \sum_{j \in J} tc_i \lambda_j y_{ij} + \sum_{i \in I} \sum_{j \in J} sc_i \lambda_j y_{ij} \quad (43)$$

s.t.

(2)–(3), (41).

By introducing auxiliary variables $r_i$, $i \in I$, the model can be transformed to the new model (38), (2)–(3), (39)–(41). ∎

**Remark 1.** The model (43), (2)–(3), (41) is also considered by Castillo et al. (2009), for which a Lagrangian-relaxation heuristic was developed. The above proposition shows that this model can be convexified and solved exactly and efficiently using existing MISOCP algorithms (see Section 6.1). A similar model is also presented in Hoseinpour and Ahmadi-Javid (2016) for designing a service system with interruption risks and solved by Lagrangian relaxation. Their problem can also be cast as an MISOCP because there is a square-root function of the decision variables $y_{ij}$ as a term in the objective function, looks like $\sqrt{\sum_{j \in J} a_{ij} y_{ij} + \left( \sum_{j \in J} b_{ij} y_{ij} \right)^2}$ with $a_{ij}, b_{ij} \geq 0$, which is similarly SOC-representable as follows:

$$\sum_{j \in J} a_{ij} y_{ij}^2 + \left( \sum_{j \in J} b_{ij} y_{ij} \right)^2 \leq r_i^2 \qquad i \in I.$$

**Remark 2.** The simplified model proposed in Proposition 3 cannot be obtained for the slightly-extended case of $var(S_i) = k/\mu_i^2$ with $k \neq 1$. Indeed, in this case, the optimal value of service time $\mu_i$ is among the real solutions of the following quartic equation (also known as the equation of the fourth degree):

$$2c_i \mu_i^4 - 4c_i \left( \sum_{j \in J} \lambda_j y_{ij} \right) \mu_i^3 + \left( 2c_i \left( \sum_{j \in J} \lambda_j y_{ij} \right)^2 - 2cw_i \sum_{j \in J} \lambda_j y_{ij} \right) \mu_i^2 + 2cw_i (1-k) \left( \sum_{j \in J} \lambda_j y_{ij} \right)^2 \mu_i - cw_i (1-k) \left( \sum_{j \in J} \lambda_j y_{ij} \right)^3 = 0.$$

Although, all the four (real or and complex) solutions of this equation can be determined by closed-form algebraic expressions, the resulting formulas are very complex and seem not useful to simplify the structure of the corresponding optimization model. One can see for $k = 1$, the last two terms disappear and a quadratic equation is obtained.

## 5 Solution algorithm

In the last decade, developing efficient algorithms have extensively been considered for solving MISOCPs. There are two classes of solution algorithms for MISOCPs. In both classes, the integrality of the problem is tackled using the generic Branch-and-Bound (B&B) method. However, in the first class, named as SOCP-based Branch-and-Bound (SOCP-B&B), the continuous relaxation of the model, which is an SOCP, is solved at each node. This is an analog of the typical B&B method applied for solving a mixed-integer linear program (for more details refer to Stubbs and Mehrotra (1999) and Leyffer (2001)). In the second class, named as Outer-Approximation-based Branch-and-Bound (OA-B&B), linear outer-approximations of MISOCPs are used to solve them. The early work in this class goes back to Quesada and Grossmann (1992), in which the classical outer approximation algorithm for mixed-integer nonlinear program (Duran & Grossmann, 1986; Fletcher & Leyffer, 1994) was extended to develop an OA-B&B algorithm. Bonami et al. (2008) generalized and applied this method for solving mixed-integer convex programs. Recently, Drewes (2009), and Drewes and Ulbrich (2012) modified the algorithm proposed by Bonami et al. (2008) for solving MISOCPs. Vielma et al. (2008) considered another method for buildings outer-approximations based on the results of Ben-Tal and Nemirovskii (2001), which generate outer polyhedral relaxations of SOC constraints. In this method, unlike the earlier one, the outer approximation constraints are not (sub)gradient-based. Krokhmal and Soberanis (2010) generalize this approach for integer $p$-order cone optimization. Vielma et al. (2017) strengthened the method presented by Vielma et al. (2008).

Our proposed models in Sections 3 and 4 can be solved by CPLEX. This solver with version $\geq 11.0$ is equipped with both SOCP-B&B and OA-B&B algorithms, as well as an option to choose the best algorithm heuristically. Numerical results show that OA-B&B algorithm seems to be more efficient than SOCP-B&B because of great capability of CPLEX in solving mixed-integer linear programs (Bonami & Tramontani, 2015).

The efficiency of the algorithms, especially for solving the general model presented in Section 3, can be improved by adding some valid inequalities to the formulation and inserting improving cuts at the root node of the B&B tree. In the rest of this section, Section 5.1 proposes a set of SOC-representable valid inequalities for the proposed MISOCPs. Then, Section 5.2 presents a cut-generation procedure to generate lift-and-project cuts at the root node to enhance the formulations.

## 5.1 Valid inequalities

Let us define the set of all feasible solutions of our problem as

$F := \{(\boldsymbol{x}, \boldsymbol{y}, \boldsymbol{\mu}, \boldsymbol{\rho}, \boldsymbol{\sigma}, \boldsymbol{u}, \boldsymbol{t}) \in \{0,1\}^n \times \{0,1\}^{m \times n} \times \mathbb{R}^n \times \mathbb{R}^n \times \mathbb{R}^{m \times n} \times \mathbb{R}^n : (2)-(5), (8), (10)-(14)\}$,

where $n = |I|, m = |J|$, and each bold letter represents a vector that includes all corresponding variables denoted by the same letter. An inequality $g(\boldsymbol{z}) \leq 0$, with $\boldsymbol{z} = (\boldsymbol{x}^T, vec(\boldsymbol{y})^T, \boldsymbol{\mu}^T, \boldsymbol{\rho}^T, \boldsymbol{\sigma}^T, vec(\boldsymbol{u})^T, \boldsymbol{t}^T)^T$, is called valid for $F$ if it is satisfied by every point in the set $F$. Theorem 2 below presents some SOC-representable valid inequalities for $F$, which can help us to tighten the feasible solution of the problem. Adding these valid inequalities may help CPLEX to reach smaller optimality gaps (as shown in our computational results in Section 6).

**Theorem 2.** (SOC-representable valid inequalities for the feasible-solution set $F$).

For $i \in I$, define $\beta_i = \min\limits_{m'_i \leq \mu_i \leq M_i} \mu_i^2 \times v_i(\mu_i)$, and $\beta'_i = \max\limits_{m'_i \leq \mu_i \leq M_i} \mu_i^2 \times v_i(\mu_i)$. Then, the inequalities

$$\left\| \begin{matrix} \alpha_i \rho_i \\ 1 - \rho_i - t_i \end{matrix} \right\| \leq 1 - \rho_i + t_i \tag{44}$$

$$\left\| \begin{matrix} \alpha'_i \sum_{j \in J} \lambda_j u_{ij} \\ 1 - \rho_i - t_i \end{matrix} \right\| \leq 1 - \rho_i + t_i \tag{45}$$

$$\sum_{j \in J} \lambda_j u_{ij} \leq \alpha''_i, \tag{46}$$

are valid for the feasible set $F$ where $\alpha_i^2 = 2(1 + \beta_i)$, $\alpha'^2_i = 2(1 + \beta'_i)/\beta'_i$, and $\alpha''^2_i = \beta'_i$.

**Proof.** From (16) and (19), one can obtain $\sum_{j \in J} \lambda_j u_{ij} = \sigma_i \sum_{j \in J} \lambda_j y_{ij} = \sigma_i \mu_i \rho_i$. Hence, from (18) it follows that

$$\left( \sum_{j \in J} \lambda_j u_{ij} \right)^2 \leq \beta'_i \rho_i^2 \tag{47}$$

$$\left( \sum_{j \in J} \lambda_j u_{ij} \right)^2 \leq \beta_i \rho_i^2. \tag{48}$$

From (48) and (18) one can get (44). Similarly, from (47) and (18) the inequality (45) can be inferred. The inequality (46) also follows from (19) and (4), which completes the proof. ∎

**Remark 3.** Note that for our case, where $v_i(\mu_i)$ is defined in (26) and is a continuous function of $\mu_i$, there always exist some values of $\beta_i, \beta'_i$ such that $0 < \beta_i \leq \beta'_i$. Moreover, even for the case where $\beta'_i \to +\infty$ and (or) $\beta_i \to 0$ the inequalities are not trivial, that is, they cannot be obtained from the existing inequalities.

Our numerical study shows that adding theses inequalities to the model are sometimes useful. In fact, the valid inequalities given in Theorem 2 tighten the continuous relaxation of $F$, denoted by $F_{\text{RL}}$ and defined as

$F_{\text{RL}} = \{(x, y, \mu, \rho, \sigma, u, t) \in [0,1]^n \times [0,1]^{m \times n} \times \mathbb{R}^n \times \mathbb{R}^n \times \mathbb{R}^{m \times n} \times \mathbb{R}^n : (2) - (5), (8), (10) - (14)\}$.

## 5.2 Enhancing lift-and-project cuts

Enhancing an MISOCP with appropriate structural cuts can considerably improve the performance of the B&B algorithm. In recent years, various kinds of cuts have been developed for MISOCPs inspiring from the ones developed for mixed-integer linear programming. Stubbs and Mehrotra (1999) develop the lift-and-project cuts of Balas et al. (1993) for 0-1 linear program to 0-1 convex problems. Cezik and Iyengar (2005) extend the lift-and-project cuts to 0-1 self-dual cone programs and use the Chvatal-Gomory procedure (Gomory, 1958) for generating linear cuts for general mixed integer self-dual cone programs.

Such cuts are specialized for 0-1 MISOCPs by Drewes (2009) and Drewes and Ulbrich (2012). Atamturk and Narayanan (2010) study a generic lifting procedure for mixed integer cone programs. Atamturk and Narayanan (2011) extend the rounding cuts of Nemhauser and Wolsey (1990) to MISOCPs. Inspiring by the work of Balas (1971, 1979) and Balas et al. (1993) for mixed-integer linear programming, a considerable effort has also been made toward generating similar kinds of cuts for MISOCPs, see, for examples, Dadush (2011), Bonami (2011), Kılınç-Karzan and Yıldız (2015), and Modaresi et al. (2016), and references therein.

In the current work, the focus is on generating lift-and-project cuts for our MISOCPs. These cuts were developed by Stubbs and Mehrotra (1999) for mixed 0-1 convex programs and further studied by Kilinç et al., (2010), and Bonami (2011). Generating the lift-and-project cuts require solving nonlinear convex programs in lifted spaces, which are generally not guaranteed to be solved in polynomial time. Recently, Drewes (2009) has applied the algorithm by Stubbs and Mehrotra (1999) to 0-1 MISOCPs. Because of the polynomial-time algorithms available for solving SOCPs even in lifted spaces, the cut-generation procedure can be effective in practice. Drewes (2009) showed that generating these cuts decreased the number of created B&B nodes for solving a benchmark set of 0-1 MISOCP. To the best of our knowledge, even though SOCPs can be solved efficiently by the CPLEX solver, these lift-and-project cuts have not been applied for solving 0-1 MISOCPs yet by this solver.

Moreover, note that most of the cuts developed for mixed-integer convex program in general and MISOCPs in particular seem neither specifically for 0-1 MISOCPs nor our analysis indicates they are not helpful for our problem, e.g., the ones developed by Cezik and Iyengar (2005) or Atamturk and Narayanan (2010), resulting in some always-valid inequalities.

We now explain how the procedure of Stubbs and Mehrotra (1999) can be applied to our problem for generating lift-and-project cuts. Define $conv(F)$ as the convex hull of the feasible-solution set $F$, i.e., the smallest convex set that contains $F$. A cut with respect to a point $z \in F_{RL}$ with $z \notin conv(F)$ is defined as an inequality that is valid for any point in the convex hull of the feasible set, but it is violated by $z$. In the following, we aim to dynamically construct a description of $conv(F)$ by generating cuts for the points of $F_{RL}$ that are not in $conv(F)$. Assume that $z_k$ indicates the $k$th binary variable among *all* binary variables of our problem, where $k \in K = \{1, \cdots, n + (n \times m)\}$. Define the set $F^B$ with $B \subseteq K$ as *a subset of $F_{RL}$* where the variables $z_k$ with $k \in B$ are binary-valued, i.e.,

$$F^B = F_{RL} \cap \{z_k \in \{0,1\}, k \in B\}.$$

One can see that $F^\emptyset = F_{RL}$ and $F^K = F_{RL}$. Using a lifting procedure, one can find a higher-dimensional description of $conv(F^B)$ as follows:

$$C_B(F_{RL}) = \begin{Bmatrix} (\mathbf{z}, \mathbf{z}^{k0}, \mathbf{z}^{k1}, \lambda_{k0}, \lambda_{k1}) \in \mathbb{R}^{|K|} \times \mathbb{R}^{|K|} \times \mathbb{R}^{|K|} \times \mathbb{R} \times \mathbb{R}: \ k \in B \\ \mathbf{z} = \lambda_{k0}\mathbf{z}^{k0} + \lambda_{k1}\mathbf{z}^{k1} \\ \lambda_{k0} + \lambda_{k1} = 1, \ \lambda_{k0}, \lambda_{k1} \geq 0 \\ \mathbf{z}^{k0} \in F_{RL}: (\mathbf{z}^{k0})_k = 0, \ \mathbf{z}^{k1} \in F_{RL}: (\mathbf{z}^{k1})_k = 1 \end{Bmatrix},$$

where symbol $(\mathbf{w})_k$ denotes the $k$th element of the vector $\mathbf{w}$. Denote $P_B(F_{RL})$ as the projection of $C_B(F_{RL})$ onto $\mathbf{z}$:

$$P_B(F_{RL}) = \{\mathbf{z}: (\mathbf{z}, \mathbf{z}^{k0}, \mathbf{z}^{k1}, \lambda_{k0}, \lambda_{k1}) \in C_B(F_{RL})\},$$

which is tighter than $F_{RL}$. Theorem 3 in the following generates cuts using the above lift-and-project procedure.

**Theorem 3.** (Lift-and-project cuts for 0-1 MISOCP)

Let $\bar{\mathbf{z}}$ be a solution in $F_{RL}$ and $\hat{\mathbf{z}}$ be the optimal solution of the following minimum distance problem

$$\min_{\mathbf{z} \in P_B(F_{RL})} \|\mathbf{z} - \bar{\mathbf{z}}\|_2 \tag{49}$$

for some $B \subseteq K$. Then, $(\hat{\mathbf{z}} - \bar{\mathbf{z}})^T \mathbf{z} \geq \hat{\mathbf{z}}^T(\hat{\mathbf{z}} - \bar{\mathbf{z}})$ is a valid linear inequality for $P_B(F_{RL})$ that cuts off $\bar{\mathbf{z}}$ whenever $\bar{\mathbf{z}} \notin P_B(F_{RL})$.

**Proof.** See Stubbs and Mehrotra (1999). ∎

The bilinear equality constraint $\mathbf{z} = \lambda_{k0}\mathbf{z}^{k0} + \lambda_{k1}\mathbf{z}^{k1}$ described in the set $C_B(F_{RL})$ is in a nonconvex format, and therefore the minimum distance problem (49) described in Theorem 3 is a nonconvex model. By introducing new variables $\mathbf{z}'^{k0} = \lambda_{k0}\mathbf{z}^{k0}, \mathbf{z}'^{k1} = \lambda_{k1}\mathbf{z}^{k1}$, the optimal solution of Problem (49) becomes equal to the optimal solution of the following problem:

$$\min_{\mathbf{z} \in P'_B(F_{RL})} \|\mathbf{z} - \bar{\mathbf{z}}\|_2, B \subseteq K \tag{50}$$

where $P'_B(F_{RL}) := \left\{\mathbf{z}: \left(\mathbf{z}, \frac{\mathbf{z}'^{k0}}{\lambda_{k0}}, \frac{\mathbf{z}'^{k1}}{\lambda_{k1}}, \lambda_{k0}, \lambda_{k1}\right) \in C_B(F_{RL})\right\}$. Proposition 4 shows that the problem (50) can be convexified as an SOCP.

**Proposition 4.** The optimization problem (50) can be reformulated as an SOCP when the service times are in the form of (24).

**Proof:** See Appendix B. ∎

The valid inequalities given in Theorem 2, and the cuts in Theorem 3, which can be generated using Proposition 1, will be applied to tighten the continuous-relaxation of the feasible set of the proposed MISOCP (i.e., $F_{RL}$) to speed up our solution algorithm. The impact of these enhancements is evaluated in the next section.

# 6 Computational results

We implemented our models in IBM ILOG OPL (Optimization Programming Language) software used the solver CPLEX 12.2, on a PC with a dual-core 2.9 GHz processor and 30GB RAM, operating Windows 7, 64-bit. To solve MISOCPS by CPLEX, the option OA-B&B was selected, since it is often more efficient than SOCP-B&B. Our B&B algorithm is enhanced by the valid inequalities proposed in Section 5.1 and the cut-generation procedure discussed in Section 5.2.

## 6.1 Performance of the algorithm

In this section, we first explain how the test problems are generated, and then analyze the performance of the algorithm that is used to solve our proposed MISOCPs under different settings.

*6.1.1 Generating test instances*

We adapted a set of well-known benchmark test problems, initially used by Holmberg et al. (1999) for a capacitated facility location problem. Similarly, the medium-size test problems in this benchmark set were adapted and used by Elhedhli (2006) for an SSD problem. The parameters of the adapted test problems such as fixed costs for establishing facilities, $ec_i$, demand rates, $\lambda_j$, and transportation costs, $tc_{ij}$, were generated according to the values of the respective parameters in the original instances. The service costs, $sc_{ij}$, and waiting costs per time unit, $wc_i$, were independently and uniformly sampled from the intervals [1, 5] and [50, 300], respectively.

*6.1.2 Performance of the solution algorithm under different settings to solve the general formulation*

Table 1 summarizes our computational results of 35 test problems using the B&B algorithm applied to solve the general MISOCP proposed in Section 3 under six different settings that indicate whether and how valid inequalities given in Theorem 2 and the cuts in Theorem 3 are incorporated. The algorithm is used for solving MISOCP (2)–(6), (8)–(14) considering (26) with $L = 2$ which covers the majority of important service-time distributions. In other words, we set $var(S_i) = \alpha_{i,0} + \frac{\alpha_{i,1}}{\mu_i^2} + \frac{\alpha_{i,2}}{\mu_i^4}$ where $\alpha_{i,l} = \delta_{i,l}^2$. $\alpha_{i,0}$ is a random from interval $(0, 1/\beta)$ where $\beta$ is square root of sum of demand rates. $\alpha_{i,i}, \alpha_{i,i}$ are randomly and independently drawn from $(0,2)$. Table 2 lists the best setting in terms of the optimality gap and number of B&B nodes.

The computational results are reported in Table 1. From this table one can see that adding the Valid Inequities (VIs) to the B&B algorithm increases the number of cuts generated by the CPLEX and reduces the possibility of occurring out-of-memory status. It can help to improve the optimality gaps in a number of instances. Adding the cuts with $|B| = 1$ completely solves the out-of-memory problem and decreases the number of nodes in most instances. It could also yield better optimality gaps in some cases. Under the other settings, the out-of-memory status frequently happens and the optimality gap can be improved only in some cases.

According to Table 2, the frequency of getting the best optimality gap under each algorithm setting is given by 8/35, 6/35, 11/35, 2/35, 7/35, and 1/35 (in order of appearance of settings in Table 1 from left to right). This clearly shows that considering either VIs or cuts (with $|B| = 1$ or $|B| = 2$) can possibly help to reduce optimality gaps. In fact, having both options is not helpful except for only 3 instances. These observations lead us to conclude that the best setting may be to consider cuts with $|B| = 1$ because it can provide the best optimality gap for 11 instances and completely prevent an out-of-memory situation for all instances. In fact, to find the best optimality gaps, it is mostly sufficient to solve the problem only for the four settings in Table 1 that do not consider both VIs and cuts simultaneously.

As can be seen in Table 2, for many cases the settings with the best optimality gap and minimum number of nodes are different. This implicitly shows that the run times required to solve the subproblems at the B&B nodes are significantly variable. In fact, there are several instances that the algorithm cannot finish solving the subproblem at the root node after three hours. Table 2 also reports the total cost and its component percentages for the best solution found for each instance. All these percentages are significant, which shows that our test problems are fairly generated.

*6.1.3 Evaluation of the simplified models for the affine and constant cases*

In Section 4, it is shown that the general model proposed in Section 3 is simplified by reducing the number of variables and constraints for two special cases of service-time distributions. Here, we assess whether these alternative formulations work better or not. To solve the model, the B&B enhanced by the cuts in Theorem 2 with $|B| = 1$ is used, as it has been considered one of the best algorithm settings in the previous subsection.

Table 3 compares the computational results obtained within the three-hour time limit for the simplified MISOCP (2)–(6), (12), (14), (34)–(35), presented in Proposition 2 for the affine case and the general model (2)–(6), (8)–(14) for six test problems with different sizes, i.e., p47, p48, p49, p54, p55, and p64. From this table, it can be observed that the simplified formulation always performs about the same or better.

Table 4 similarly evaluates the model (2)–(3), (38)–(41) developed for the constant case with $m_i = 0$, $M_i \to +\infty$. The results show that this model significantly outperforms the general model and can find the exact or near optimal solutions, while the general model provides highly suboptimal solutions.

**Table 1.** Computational results of the algorithm under six different settings for solving the general model proposed in Section 3 within a three-hour run-time limit.

| Test Problem | Number of SFs | Number of Customers | Basic B&B | | | B&B + VIs | | | B&B + cuts for $|B|=1$ | | | B&B + VIs + cuts for $|B|=1$ | | | B&B + cuts for $|B|=2$ | | | B&B + VIs + cuts for $|B|=2$ | | |
|---|---|---|---|---|---|---|---|---|---|---|---|---|---|---|---|---|---|---|---|---|
| | | | Number of nodes | Number of cuts | Gap (%) | Number of nodes | Number of cuts | Gap (%) | Number of nodes | Number of cuts | Gap (%) | Number of nodes | Number of cuts | Gap (%) | Number of nodes | Number of cuts | Gap (%) | Number of nodes | Number of cuts | Gap (%) |
| p1 | 10 | 50 | <273k | 5 | 10.14 | <422k | 1 | 10.82 | <75k | 127 | 0.16 | <679k | 221 | 4.44 | <50k | 125 | 0.11 | - | - | - |
| p2 | 10 | 50 | <241k | 0 | 13.16 | - | - | - | <187k | 181 | 8.59 | <151k | 259 | 17.67 | <171k | 247 | 8.86 | <195k | 331 | 16.99 |
| p3 | 10 | 50 | <291k | 2 | 11.14 | <254k | 6 | 13.26 | <134k | 273 | 8.79 | - | - | - | - | - | - | - | - | - |
| p4 | 10 | 50 | <253k | 6 | 9.58 | <1156k | 28 | 10.48 | <87k | 243 | 9.80 | <127k | 289 | 12.62 | <85k | 253 | 9.65 | <188k | 302 | 12.46 |
| p13 | 20 | 50 | - | - | - | <302k | 24 | 23.40 | <127k | 314 | 17.35 | - | - | - | <102k | 376 | 17.66 | - | - | - |
| p14 | 20 | 50 | <168k | 4 | 20.52 | <195k | 27 | 25.39 | <127k | 371 | 9.01 | <179k | 413 | 9.30 | <122k | 372 | 9.05 | <144k | 423 | 9.43 |
| p15 | 20 | 50 | - | - | - | - | - | - | <51k | 357 | 11.28 | <148k | 452 | 12.23 | <49k | 422 | 12.04 | <105k | 408 | 6.59 |
| p16 | 20 | 50 | - | - | - | <218k | 160 | 18.86 | <75k | 482 | 14.53 | <105k | 642 | 16.99 | <77k | 460 | 14.30 | <89k | 655 | 17.66 |
| p25 | 30 | 150 | - | - | - | - | - | - | 78 | 494 | 17.04 | 612 | 524 | 17.14 | - | - | - | - | - | - |
| p26 | 30 | 150 | - | - | - | 121 | 341 | 7.77 | 660 | 416 | 11.98 | - | - | - | <51k | 778 | 14.73 | - | - | - |
| p27 | 30 | 150 | <336k | 88 | 8.47 | - | - | - | <49k | 822 | 9.36 | - | - | - | <52k | 919 | 8.39 | - | - | - |
| p28 | 30 | 150 | <574k | 139 | 7.42 | - | - | - | <29k | 506 | 5.81 | 328 | 833 | 10.71 | <36 | 486 | 4.87 | 136 | 574 | 9.11 |
| p41 | 10 | 90 | <59k | 184 | 11.40 | <63k | 274 | 11.54 | <63k | 228 | 9.51 | <64k | 237 | 11.82 | <69k | 238 | 9.66 | - | - | - |
| p42 | 20 | 80 | - | - | - | - | - | - | <6k | 735 | 8.73 | <10k | 852 | 11.81 | <9k | 662 | 8.79 | <17k | 761 | 10.65 |
| p43 | 30 | 70 | <26k | 250 | 7.81 | <305k | 18 | 12.15 | <3k | 407 | 11.14 | <13k | 484 | 12.00 | <3k | 407 | 11.14 | <13k | 484 | 12.00 |
| p44 | 10 | 90 | <80k | 271 | 13.48 | <67k | 232 | 14.96 | <108k | 522 | 13.61 | <186k | 459 | 14.82 | - | - | - | <177k | 474 | 15.35 |
| p45 | 20 | 80 | <188k | 3 | 12.00 | <207k | 3 | 12.69 | <9k | 480 | 6.58 | <9k | 421 | 9.71 | <10k | 473 | 7.72 | <24k | 465 | 9.05 |
| p46 | 30 | 70 | <258k | 2 | 9.33 | <26k | 485 | 5.70 | <59k | 156 | 10.55 | <25k | 417 | 8.76 | <80k | 217 | 10.75 | - | - | - |
| p47 | 10 | 90 | <243k | 0 | 19.53 | - | - | - | <84k | 274 | 14.17 | <77k | 333 | 17.21 | <79k | 279 | 14.01 | <77k | 340 | 16.86 |
| p48 | 20 | 80 | - | - | - | <17k | 508 | 8.97 | <79k | 727 | 16.00 | <105k | 981 | 15.98 | <82k | 805 | 15.28 | <78k | 987 | 16.66 |
| p49 | 30 | 70 | - | - | - | <362k | 8 | 13.62 | <6k | 611 | 11.26 | <95k | 738 | 12.62 | <7k | 427 | 12.64 | <112k | 1011 | 12.56 |
| p50 | 10 | 100 | <58k | 103 | 10.89 | <49k | 145 | 11.99 | <49k | 197 | 4.30 | <48k | 266 | 6.29 | <49k | 173 | 4.59 | <50k | 194 | 6.43 |
| p51 | 20 | 100 | <6k | 329 | 5.50 | <7k | 260 | 7.95 | <7k | 246 | 1.34 | <7k | 413 | 3.93 | <8k | 248 | 1.33 | <11k | 500 | 3.70 |
| p52 | 10 | 100 | <60k | 294 | 13.87 | <65k | 231 | 15.47 | <57k | 159 | 3.56 | <76k | 200 | 4.31 | <117k | 166 | 5.29 | <172k | 217 | 4.29 |
| p53 | 20 | 100 | <7k | 339 | 6.96 | <9k | 218 | 9.21 | <73k | 683 | 10.75 | <8k | 491 | 14.39 | <59k | 733 | 10.80 | <8k | 486 | 14.48 |
| p54 | 10 | 100 | <76k | 317 | 15.39 | <60k | 282 | 18.25 | <73k | 513 | 13.66 | <151k | 494 | 16.04 | <61k | 479 | 13.61 | <64k | 344 | 16.18 |
| p55 | 20 | 100 | <245k | 6 | 13.31 | - | - | - | <4k | 475 | 15.55 | <84k | 767 | 14.78 | <4k | 599 | 15.04 | - | - | - |
| p64 | 30 | 200 | 0 | 424 | 30.21 | 186 | 793 | 12.68 | 0 | 731 | 30.19 | 6 | 870 | 35.83 | 0 | 683 | 11.12 | 0 | 831 | 44.84 |
| p65 | 30 | 200 | 0 | 482 | 20.29 | 0 | 617 | 52.42 | 0 | 857 | 28.75 | 10 | 882 | 6.66 | 0 | 878 | 22.18 | 0 | 894 | 36.57 |
| p66 | 30 | 200 | - | - | - | - | - | - | <34k | 1983 | 18.51 | <42k | 2345 | 15.54 | <32k | 2566 | 14.18 | <45k | 2503 | 14.36 |
| p67 | 30 | 200 | 0 | 475 | 11.76 | 0 | 463 | 28.54 | 20 | 1383 | 35.92 | 0 | 1117 | 37.70 | 0 | 890 | 36.30 | 0 | 910 | 71.81 |
| p68 | 30 | 200 | 0 | 478 | 14.07 | 0 | 633 | 29.70 | 0 | 532 | 16.60 | 0 | 626 | 30.59 | 0 | 529 | 34.23 | 0 | 616 | 18.87 |
| p69 | 30 | 200 | 0 | 496 | 15.48 | 53 | 740 | 11.58 | 0 | 935 | 32.85 | 0 | 882 | 24.13 | 0 | 885 | 22.39 | 0 | 899 | 48.28 |
| p70 | 30 | 200 | - | - | - | - | - | - | <50k | 1557 | 11.94 | <66k | 2380 | 9.17 | <34k | 1448 | 12.65 | <58k | 2045 | 10.44 |
| p71 | 30 | 200 | 73 | 531 | 7.8 | 0 | 381 | 50.3 | 8 | 394 | 10.6 | 24 | 632 | 10.3 | <46k | 975 | 9.9 | 14 | 551 | 11.9 |

- Out-of-memory before reaching the time limit

**Table 2. The total cost and it component percentages for the best solution, and the best algorithm setting in terms of runt time and number of nodes for the results reported in Table 1.**

| Test Problem | Number of SFs | Number of Customers | Total Cost | Serving Cost (%) | Waiting Cost (%) | Traveling Cost (%) | Opening Cost (%) | Best Gap (%) | Algorithm setting resulting in best gap | Algorithm setting resulting in minimum number of nodes** |
|---|---|---|---|---|---|---|---|---|---|---|
| p1 | 10 | 50 | 2703 | 28 | 10 | 29 | 33 | 0.11 | B&B + cuts for $|B|=2$ | B&B + cuts for $|B|=2$ |
| p2 | 10 | 50 | 1958 | 19 | 18 | 32 | 31 | 8.59 | B&B + cuts for $|B|=1$ | B&B + VIs + cuts for $|B|=1$ |
| p3 | 10 | 50 | 2312 | 16 | 12 | 29 | 43 | 8.79 | B&B + cuts for $|B|=1$ | B&B + cuts for $|B|=1$ |
| p4 | 10 | 50 | 2966 | 22 | 9 | 22 | 47 | 9.65 | B&B + cuts for $|B|=1$ | B&B + cuts for $|B|=2$ |
| p13 | 20 | 50 | 2324 | 31 | 15 | 18 | 37 | 17.35 | B&B + cuts for $|B|=1$ | B&B + cuts for $|B|=2$ |
| p14 | 20 | 50 | 2348 | 41 | 19 | 15 | 26 | 9.01 | B&B + cuts for $|B|=1$ | B&B + cuts for $|B|=2$ |
| p15 | 20 | 50 | 2970 | 37 | 18 | 11 | 34 | 6.59 | B&B + VIs + cuts for $|B|=2$ | B&B + cuts for $|B|=2$ |
| p16 | 20 | 50 | 3148 | 30 | 14 | 11 | 44 | *13.79 | Basic B&B | B&B + cuts for $|B|=1$ |
| p25 | 30 | 150 | 4280 | 28 | 11 | 30 | 31 | *7.47 | B&B + VIs | B&B + cuts for $|B|=1$ |
| p26 | 30 | 150 | 3042 | 34 | 12 | 34 | 20 | 7.77 | B&B + VIs | B&B + VIs |
| p27 | 30 | 150 | 4151 | 28 | 12 | 35 | 24 | *5.77 | B&B + VIs | B&B + cuts for $|B|=1$ |
| p28 | 30 | 150 | 4472 | 26 | 11 | 32 | 31 | 4.87 | B&B + cuts for $|B|=2$ | B&B + VIs + cuts for $|B|=2$ |
| p41 | 10 | 90 | 3246 | 27 | 13 | 35 | 26 | 9.51 | B&B + cuts for $|B|=1$ | Basic B&B |
| p42 | 20 | 80 | 2577 | 22 | 12 | 41 | 25 | 8.73 | B&B + cuts for $|B|=1$ | B&B + cuts for $|B|=1$ |
| p43 | 30 | 70 | 1825 | 15 | 12 | 40 | 33 | 7.81 | Basic B&B | B&B + cuts for $|B|=1$ |
| p44 | 10 | 90 | 2860 | 30 | 14 | 34 | 22 | 13.48 | Basic B&B | B&B + cuts for $|B|=1$ |
| p45 | 20 | 80 | 2218 | 24 | 7 | 37 | 33 | 6.58 | B&B + cuts for $|B|=1$ | B&B + cuts for $|B|=1$ |
| p46 | 30 | 70 | 1845 | 11 | 14 | 32 | 42 | 5.70 | B&B + VIs | B&B + cuts for $|B|=1$ |
| p47 | 10 | 90 | 3182 | 29 | 10 | 31 | 29 | 14.01 | B&B + cuts for $|B|=2$ | B&B + VIs + cuts for $|B|=2$ |
| p48 | 20 | 80 | 2050 | 17 | 14 | 35 | 33 | 8.97 | B&B + VIs | B&B + VIs |
| p49 | 30 | 70 | 2122 | 30 | 11 | 28 | 31 | 11.26 | B&B + cuts for $|B|=1$ | B&B + cuts for $|B|=1$ |
| p50 | 10 | 100 | 3597 | 22 | 10 | 45 | 23 | 4.30 | B&B + cuts for $|B|=1$ | B&B + VIs + cuts for $|B|=1$ |
| p51 | 20 | 100 | 2785 | 15 | 5 | 51 | 28 | 1.33 | B&B + cuts for $|B|=2$ | Basic B&B |
| p52 | 10 | 100 | 4428 | 27 | 7 | 37 | 29 | 3.56 | B&B + cuts for $|B|=1$ | B&B + cuts for $|B|=1$ |
| p53 | 20 | 100 | 2562 | 21 | 11 | 45 | 24 | 6.96 | Basic B&B | B&B + cuts for $|B|=2$ |
| p54 | 10 | 100 | 3388 | 27 | 9 | 37 | 27 | 13.61 | B&B + cuts for $|B|=2$ | B&B + VIs |
| p55 | 20 | 100 | 2146 | 15 | 11 | 41 | 32 | 13.31 | Basic B&B | B&B + cuts for $|B|=2$ |
| p64 | 30 | 200 | 5637 | 33 | 6 | 43 | 18 | 11.12 | B&B + cuts for $|B|=2$ | B&B + VIs + cuts for $|B|=1$ |
| p65 | 30 | 200 | 6159 | 32 | 7 | 35 | 26 | 6.66 | B&B + VIs + cuts for $|B|=1$ | Basic B&B |
| p66 | 30 | 200 | 7898 | 24 | 10 | 29 | 38 | 14.18 | B&B + cuts for $|B|=2$ | B&B + cuts for $|B|=1$ |
| p67 | 30 | 200 | 5035 | 18 | 9 | 44 | 30 | 11.76 | Basic B&B | B&B + cuts for $|B|=1$ |
| p68 | 30 | 200 | 4570 | 20 | 8 | 50 | 22 | 14.07 | Basic B&B | Basic B&B |
| p69 | 30 | 200 | 5189 | 19 | 10 | 40 | 31 | 11.58 | B&B + VIs | Basic B&B |
| p70 | 30 | 200 | 8488 | 30 | 8 | 27 | 35 | 9.17 | B&B + VIs + cuts for $|B|=1$ | B&B + cuts for $|B|=2$ |
| p71 | 30 | 200 | 4850 | 19 | 9 | 44 | 28 | 7.8 | Basic B&B | B&B + cuts for $|B|=1$ |

* Out-of-memory before reaching the time limit

** When the number of non-root nodes is zero, the setting with the best gap is reported

**Table 3.** Evaluation of the simplified model (2)–(6), (12), (14), (34)–(35) for the affine case of service-time distributions.

| Test Problem | Number of SFs | Number of Customers | Gap % | |
|---|---|---|---|---|
| | | | General model | Simplified model |
| p48 | 20 | 80 | 15.27 | 4.27 |
| p49 | 30 | 70 | 10.96 | 10.92 |
| p54 | 10 | 100 | 13.44 | 13.41 |
| p55 | 20 | 100 | 12.45 | 12.35 |
| p64 | 30 | 200 | 28.25 | 7.65 |

**Table 4.** Evaluation of the simplified model (2)–(3), (38)–(41) for the constant case of service-time distributions.

| Test Problem | Number of SFs | Number of Customers | Gap % | |
|---|---|---|---|---|
| | | | General model | Simplified model |
| p48 | 20 | 80 | 17.13 | *0.00 |
| p49 | 30 | 70 | 9.37 | *0.00 |
| p54 | 10 | 100 | 12.55 | *0.00 |
| p55 | 20 | 100 | 15.91 | *0.00 |
| p64 | 30 | 200 | 35.99 | 0.73 |

* Optimal integer solution is found.

### 6.2 A real-world example: Locating Preventive Medical Facilities in Toronto, Canada

In the last years, there has been an increased interest in preventive healthcare programs. In this example, we analyze locating a set of preventive mammography clinics in city of Toronto and design a service network to respond the needs of the targeted people. Our main goal is to show that recognizing the correct service time distributions at SFs, at least their means and variances, are very important since various variance structures can result in a different network design for the service system.

The 96 large regions of Toronto, called FSAs using the Canada Post classification, are selected for providing service. We selected 19 general service hospitals in the city as potential places for locating the preventive clinics. Euclidean travel distances between the region centroids were recorded from Statics Canada. Demographic data at the FSA level were also available from Statistics Canada. The targeted people, who require service once a year, are the females aged 50–59, 0.071% of Canadian population, refer to Aboolian et al. (2015). Participation rate is assumed to be no more than 95% based on Verter and Lapierre (2002). Other parameters of the problem are randomly generated from some logical intervals.

Figure 1 shows the location-allocation decisions graphically. As seen in this figure, three SFs are established to serve the people. Table 5 reports the mean and variance of service times at each open SF as well as number of the customers served by each one. Percentages of cost terms; i.e., service, waiting, traveling and establishing cost, in the objective function are 15.5%, 26.3%, 32.6%, and 25.6%, respectively.

**Table 5. Results for the case study**

| FSAs with open SF | Number of assigned FSAs | Mean of service time | Variance of service time |
|---|---|---|---|
| M2M | 36 | 355 | 0.0013 |
| M5A | 31 | 220 | 0.0019 |
| M6S | 29 | 300 | 0.0029 |

Increasing the coefficient of travel and waiting costs yields to the optimal solution with more open SFs in which most of the customers are assigned to SFs that are less congested. Moreover, about 90% of customers are assigned to the closest open SFs, for which constraints (7) are mostly satisfied, while they are not considered here. This shows that the user-choice and central-authority assignment policies are highly correlated.

Figure 2 shows the effect of the type of service-time distribution on the designed network structure; i.e., the optimal location-allocation decisions. As shown in the Figure 2.a, when the service-time distributions in all SFs are assumed to be exponential, the same number of facilities are decided to open, but in different locations. When all the service time distributions are assumed to be Gamma distributions, the number and location of the established SFs are significantly changed; see Figure 2.b. As a direct result of the change in location of open SFs, the allocation decisions, customer assignments, will also be changed. This shows that in order to better model the real world SSDPs, estimating the correct service time distributions is highly important.

## 7 A unified convexification approach to metrics of M/G/1 queue

The goal of this section is to show the general applicability of the method used in the paper for constructing mixed-integer convex optimization models. It provides various convexity results and discussions for M/G/1 queue metrics. Only a part of the results given in Theorem 5 below has been used in Sections 2 and 3, and the other parts are new. In the following, we first discuss the possibility of convex reformulation, and then SOCP reformulation in particular.

Assume that $\mu$ is a continuous variable (satisfying $\mu \geq \lambda$), $v(\mu)$ is a variance function, and $\lambda$ is represented by $\sum_{i=1}^{n} \lambda_i w_i$ with $\mathbf{w} = (w_1, \ldots, w_n)^T \in \{0,1\}^n$ where $\lambda_i$, $i = 1, \ldots, n$ are given positive constants. Now we investigate sufficient conditions under which the following constraints:

$$WT_T(\lambda, \mu) \leq z \tag{51}$$

$$WT_I(\lambda, \mu) \leq z \tag{52}$$

can be convexified in terms of $\mu$ and $\mathbf{w}$ (or $\mu$ and $\lambda$) where $z \geq 0$ is a constant, a variable, or any affine function of $\mu$ and $\lambda$; and where the quantities

$$WT_T(\lambda, \mu) = \lambda \, WT_I(\lambda, \mu) = \frac{\lambda^2(1 + v(\mu)\mu^2)}{2\mu(\mu - \lambda)} + \frac{\lambda}{\mu}$$

$$WT_I(\lambda, \mu) = \frac{\lambda(1 + v(\mu)\mu^2)}{2\mu(\mu - \lambda)} + \frac{1}{\mu}$$

are the *total* and *individual* expected waiting times in an M/G/1 queue system, respectively. Note that $WT_T(\lambda, \mu)$ is also the expected number of people waiting in the systems. Constraint (51) is studied in the next theorem.

**Theorem 4.** If constraint $v(\mu) \leq \sigma^2$ can be convexified, then constraint (51) can be convexified in terms of $\mu$ and $\boldsymbol{w}$ as

$$\sqrt{4\sum_{i=1}^{n} \lambda_i w_i^2 + (\rho - \mu)^2} \leq \rho + \mu$$

$$\frac{\rho^2 + (\sum_{i=1}^{n} \lambda_i u_i)^2}{2(1 - \rho)} + \rho \leq z$$

$$\sigma - (1 - w_i)M \leq u_i \leq \sigma$$

$$0 \leq u_i \leq w_i M$$

$$\mu \geq \sum_{i=1}^{n} \lambda_i w_i, \ 0 \leq \rho \leq 1, \boldsymbol{w} \in \{0,1\}^n, \boldsymbol{u} \geq \boldsymbol{0}$$

where $M$ is a sufficiently large positive constant. If $v(\mu) \leq \sigma^2$ is SOC-representable, then this reformulation is also SOC-representable. In particular, if $v(\mu) = a + b\mu^{-2}$ with $a, b \geq 0$, then constraint (51) can be convexified more efficiently (without needing the $n$ auxiliary variables $u_i$, $i = 1, \ldots, n$ and $4n$ associated constraints) as follows:

$$\sqrt{4\sum_{i=1}^{n} \lambda_i w_i^2 + (\rho - \mu)^2} \leq \rho + \mu$$

$$\frac{\rho^2(1 + b) + a(\sum_{i=1}^{n} \lambda_i w_i)^2}{2(1 - \rho)} + \rho \leq z$$

$$\mu \geq \sum_{i=1}^{n} \lambda_i w_i, \ 0 \leq \rho \leq 1, \boldsymbol{w} \in \{0,1\}^n,$$

which is also SOC-representable.

**Proof:** The proof is similar to the ones presented for Theorem 1 and Proposition 1. Recall that in the first part, we need to define the auxiliary variables $u_i = \sigma w_i$, $i = 1, \ldots, n$, and $\rho = \lambda/\mu$. ∎

Note that the above results are used in Sections 2 and 3 to solve our decision problem. The next theorem investigates the case where the individual waiting time in an M/G/1 queue system is used as a performance measure.

**Theorem 5.** If constraint $v(\mu) \leq \sigma^2$ can be convexified, then constraint (52) can be convexified in terms of $\mu$ and $\boldsymbol{w}$ as

$$\sqrt{4\sum_{i=1}^{n} \lambda_i w_i^2 + (\rho - \mu)^2} \leq \rho + \mu$$

$$\frac{\sum_{i=1}^{n} \lambda_i w_i^2}{2(\mu - \sum_{i=1}^{n} \lambda_i w_i)} + \frac{\sum_{i=1}^{n} \lambda_i u_i^2}{2(1-\rho)} + \frac{1}{2\mu} \leq z$$

$$\sigma - (1 - w_i)M \leq u_i \leq \sigma$$

$$0 \leq u_i \leq w_i M$$

$$\mu \geq \sum_{i=1}^{n} \lambda_i w_i,\ 0 \leq \rho \leq 1,\ \boldsymbol{w} \in \{0,1\}^n,\ \boldsymbol{u} \geq \boldsymbol{0}.$$

In particular, if $v(\mu) = a + b\mu^{-2}$ with $a, b \geq 0$, then constraint (52) can be convexified more efficiently in terms of $\mu$ and $\boldsymbol{w}$:

$$\sqrt{4\sum_{i=1}^{n} \lambda_i w_i^2 + (\rho - \mu)^2} \leq \rho + \mu$$

$$\frac{(1+b)\sum_{i=1}^{n} \lambda_i w_i^2}{2(\mu - \sum_{i=1}^{n} \lambda_i w_i)} + \frac{a\sum_{i=1}^{n} \lambda_i w_i^2}{2(1-\rho)} + \frac{1}{2\mu} \leq z$$

$$\mu \geq \sum_{i=1}^{n} \lambda_i w_i,\ 0 \leq \rho \leq 1,\ \boldsymbol{w} \in \{0,1\}^n,$$

which is SOC-representable.

**Proof:** To obtain the results, one needs to consider the identity

$$WT_I(\lambda, \mu) = \frac{b+1}{2(\mu - \lambda)} + \frac{a\lambda}{2(1-\rho)} + \frac{1-b}{2\mu}$$

with $\rho = \lambda/\mu$, which follows (after some work) from

$$\frac{\lambda}{2\mu(\mu-\lambda)} = \frac{1}{2(\mu-\lambda)} - \frac{1}{2\mu}.$$

Moreover, the function $x^2/y$ with $y > 0$ is convex, so its composition with affine functions is convex. The rest of the proof is similar and straightforward. ∎

**Remark 4.** It worth nothing that if $\mu$ depends on another vector of decision variables $\boldsymbol{x}$ as $\mu = h(\boldsymbol{x})$, the latter can be relaxed as $\mu \leq h(\boldsymbol{x})$ because queue waiting-time measures are decreasing in $\mu$. Hence, if the function $h$ is concave, the above formulations with $\mu \leq h(\boldsymbol{x})$ remain convex.

**Remark 5.** One should carefully consider that the results derived in Ahmadi-Javid and Hoseinpour (2017) is for a special case that $\lambda = \sum_{i=1}^{n} \lambda_i w_i$ with $\boldsymbol{w} \in \{0,1\}^n$ and $\mu$ and $\sigma^2$ is represented by $\sum_{j=1}^{m} \mu_j t_j$ and $\sum_{j=1}^{m} \sigma_j^2 t_j$ with $\boldsymbol{t} \in \{0,1\}^m$ where $\mu_j$ and $\sigma_j^2\ j = 1, \ldots, m$ are given positive constants. Indeed, for this special case based on the above general results one can obtain new convex formulations, but with less computational performance because of big-M constraints used in them.

In both of the above theorems, it is assumed that $\lambda = \sum_{i=1}^{n} \lambda_i w_i$ with $\mathbf{w} \in \{0,1\}^n$, which means that $\lambda$ must be selected from a finite set of given alternatives. Unfortunately, providing sufficient conditions to guarantee that constraints (51) and (52) can be convexified in terms of $\mu$ and $\lambda$ is challenging. A few of such conditions are presented in the next theorem.

**Theorem 6.** If $v(\mu) = b\mu^{-2}$ with $b \geq 0$ and $z$ is fixed to a constant $c$, then constraint (51) can be convexified in terms of $\mu$ and $\lambda$ as

$$\frac{\lambda^2(1+b)}{2(\mu - \lambda)} + \lambda \leq c\mu, \mu \geq \lambda, \lambda \geq 0,$$

which is also SOC-representable. If $v(\mu) = a + b\mu^{-2}$ with $a \geq 0$ and $0 \leq b \leq 1$, then constraint (52) can be convexified in terms of $\mu$ and $\lambda$ as follows:

$$\frac{b+1}{2(\mu - \lambda)} + \frac{a\lambda\mu}{2(\mu - \lambda)} + \frac{1-b}{2\mu} \leq z$$

$$\mu \geq \lambda, \lambda \geq 0,$$

which is SOC-representable if $a = 0$. When $\mu$ is assumed to be fixed, both constraints (51) and (52) are convex in variable $\lambda$, but only constraint (51) is SOC-representable.

**Proof:** The proof of the first part is clear by noting that $\frac{\lambda^2}{(\mu-\lambda)}$ is the perceptive function associated with the convex function $\frac{\rho^2}{1-\rho}$. To prove the second assertion, one can see that

$$WT_I(\lambda, \mu) = \frac{b+1}{2(\mu - \lambda)} + \frac{a\lambda\mu}{2(\mu - \lambda)} + \frac{1-b}{2\mu}.$$

The convexity of the first term is clear. The second term is always convex as it is the perceptive function associated with the convex function $\frac{a\rho}{2(1-\rho)}$. The third term is convex only if $b \leq 1$. This completes the proof. ∎

**Remark 6.** When $\lambda$ is a function of a decision vector $\mathbf{y}$ with *arbitrary* domain, that is, $\lambda = g(\mathbf{y})$, the equality $\lambda = g(\mathbf{y})$ can be relaxed as $g(\mathbf{y}) \leq \lambda$ (queue waiting-time measures depend increasingly on $\lambda$). Thus, the formulations obtained above with the new constraint $g(\mathbf{y}) \leq \lambda$ remains convex for a convex function $g$.

**Remark 7.** Only in the first part of Theorem 6, we restrict our condition to the case where z is a constant, and in all the other places of the section the variable $z$ can be a constant, variable, or an affine function of variables $\mu$ and $\lambda$. In fact, when z is a constant in constraints (51) and (52), some other reformulations and simplifications can be presented, which are not listed here for the sake of shortness. Moreover, if one is interested to consider the waiting times in the queue (instead of in the system), the terms $\lambda\mu^{-1}$ and $\mu^{-1}$ must be removed from constraints (51) and (52), respectively. This removal can result in simpler reformulations.

**Remark 8.** One can see that $v(\mu) = b\mu^{-2}$ with $b \geq 0$ covers the simple representation (23) used by Weber (1983), Stidham (1992), and Fridgeirsdottir and Chiu (2005). Theorem 6 highlights that when the arrivals evolve according to a Poisson process $WT_T(\lambda, \mu)$ is not jointly convex in both $\lambda$ and $\mu$, and that $WT_I(\lambda, \mu)$ is convex in both $\lambda$ and $\mu$ only if $0 \leq b \leq 1$. This implies that for G/G/1 queues $WT_T(\lambda, \mu)$ is not generally convex in both $\lambda$ and $\mu$, while it is convex in each of $\lambda$ and $\mu$. Note that $WT_T(\lambda, \mu)$ is not generally convex in $\lambda$ for G/G/1 queues (Fridgeirsdottir & Chiu, 2005).

The above results can be used when quantities $WT_T(\lambda, \mu)$ and $WT_I(\lambda, \mu)$ appear in an optimization problem to control the overall congestion or satisfy some service levels while $\lambda$ and $\mu$ are either decision variables or depend on some other decision variables. The implementation procedure is similar to the one used in this paper.

# 8 Concluding remarks

This paper considers a balance-objective SLCLS problem with M/G/1 queues where the location and service capacity planning are integrated. It is assumed that the variance of the service time is a function of the service capacity (rate), which is the reciprocal of the mean service time. This service capacity is selected from a given (bounded or unbounded) intervals. The problem is initially formulated as a mixed-integer non-convex model, and then reformulated as a novel MISOCP, which can be solved optimally using recent advances in mixed-integer convex programming. Valid inequalities and a cut-generation procedure are developed to achieve more computational efficiency. Our numerical study on a health care location problem using some real data shows the high impact of service-time distributions on the network design of a service system.

Extending the problem under more realistic assumptions (e.g., uncertain demand rates, failure-prone servers, multi-services, elastic demand rates) is a possible direction.

The modeling method used here can be applied to other decision problems with M/G/1 queues where the service rate of each queue is a continuous decision variable and the demand rate is either known constant or a number in a finite set, which is determined by an affine function of binary/integer variables.

Some results on the case where demand rates are continuous decision variables are obtained in this paper for special cases. However, it is important to investigate if these results can be generalized. Extending our results for systems with G/G/1 queues would be very interesting but seems extremely difficult as there is no know closed-form formulas for the performance metrics of these queues and thus it may remain an open problem.

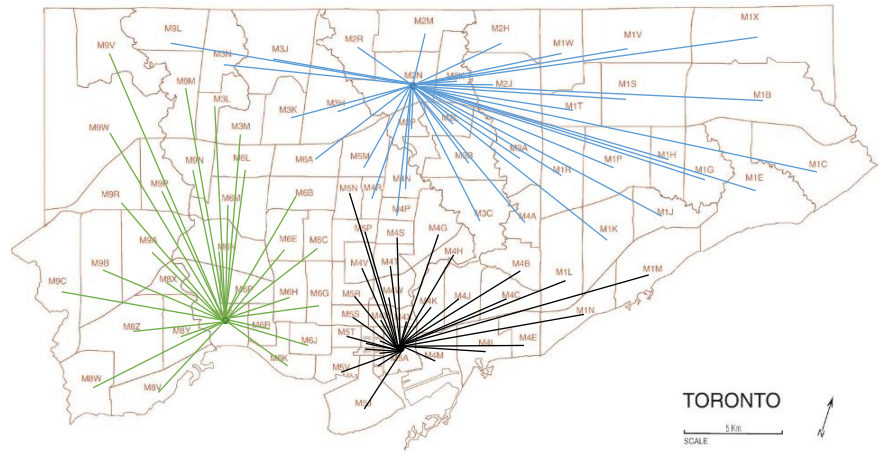

**Figure 1. Location-allocation decisions of the case study with general-form service time distributions.**

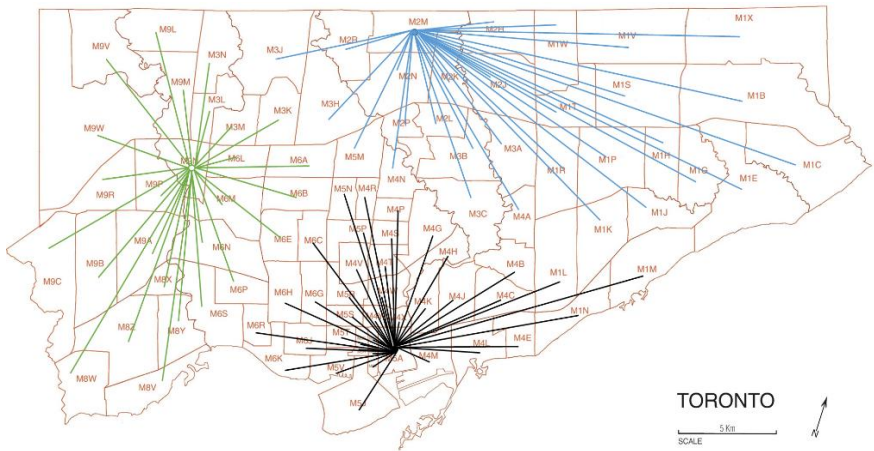

| FSAs with Open SF | Number of assigned customers | Mean of service time | Variance of service time |
|---|---|---|---|
| M2M | 30 | 333 | 0.0030 |
| M5A | 39 | 294 | 0.0034 |
| M6N | 27 | 266 | 0.0038 |

a) **Location-allocation decisions with exponential service-time distributions.**

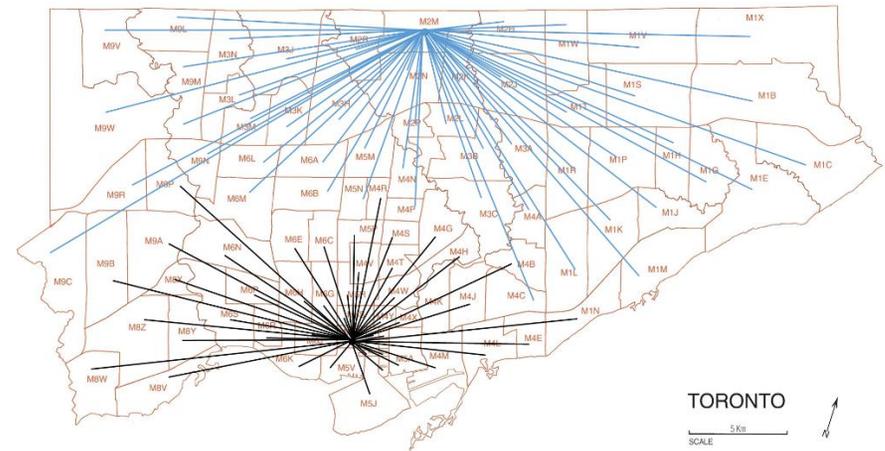

| FSAs with Open SF | Number of assigned customers | Mean of service time | Variance of service time |
|---|---|---|---|
| M2M | 50 | 577 | 0.0039 |
| M5T | 46 | 364 | 0.0061 |

b) **Location-allocation decisions with Gamma service-time distributions.**

**Figure 2. Location-allocation decisions of the case study with special service time distributions**

## Appendix A. Preliminaries

This appendix provides some preliminaries needed throughout the paper on MISOCPs. Let $\kappa \subseteq \mathbb{R}^m$ be a proper cone, i.e., a pointed and closed convex cone with a nonempty interior, which induces a partial ordering on $\mathbb{R}^m$, denoted by $\geq_\kappa$, where we have $a \geq_\kappa b \Leftrightarrow a - b \geq_\kappa 0 \Leftrightarrow a - b \in \kappa$.

The optimization problem

$$\min_{x \in \mathbb{R}^n} c^T x$$

$$Ax - b \geq_\kappa 0 \quad x \in P,$$

is called a *cone program* where $c \in \mathbb{R}^n, b \in \mathbb{R}^n, A$, is an $m \times n$ matrix, and $P$ is a polyhedron. A cone problem is an optimization model with linear objective function and generalized linear inequalities $Ax - b \geq_\kappa 0$.

A *second order cone* is a proper cone and is defined as:

$$S^m = \left\{ z \in \mathbb{R}^m : z_m \geq \sqrt{\sum_{i=1}^{m-1} z_i^2} \right\}.$$

This cone is also called a Lorentz, quadratic, or ice-cream cone. A Second-Order Cone Program (SOCP), also called a Conic Quadratic Program (CQP), is a cone program for which the cone $\kappa$ is a direct product of several second-order cones, $\kappa = S^{m_1} \times S^{m_1} \times \cdots S^{m_k}$. Therefore, an SOCP is like as

$$\min_{x \in \mathbb{R}^n} c^T x$$

$$A_i x - b_i \geq_{S^{m_i}} 0 \quad i = 1, \dots, k$$

$$x \in P.$$

By defining $[A_i, b_i] = \begin{bmatrix} E_i & e_i \\ \beta_i^T & \delta_i \end{bmatrix}$, this can simply be written as follows:

$$\min_{x \in \mathbb{R}^n} c^T x$$

$$\|E_i x - e_i\| \leq \beta_i^T x - \delta_i \quad i = 1, \dots, k$$

$$x \in P.$$

When all $\beta_i$ and $\delta_i$ are zeroes, the SOCP reduces to a linear program. If all $p_i$ are zeroes the SCOP becomes a convex quadratically constrained linear program. SOCPs are important convex optimization problems. The most popular solution methods for SOCPs are interior point methods because of their polynomial-time theoretical convergence and efficient computational performance in implementation. For more details on SOCPs and their applications see, e.g., Alizadeh and Goldfarb (2003), Alizadeh and Goldfarb (2003), Ben-Tal and Nemirovski (2001), and Boyd and Vandenberghe (2004).

Each constraint $\|E_i x - e_i\| \leq \beta_i^T x - \delta_i$ is called an SOC constraint. A set of constraints (or generally any subset of $\mathbb{R}^n$) is called SOC-representable if it can be expressed by a finite set of linear and SOC constraints. A mathematical program is called a Mixed-Integer Second-Order Cone Program (MISOCP) when some of decision variables $x_i$ $i = 1, \ldots, n$ are restricted to be integer.

## Appendix B. Proof of Proposition 4

Without loss of generality, assume that $k \leq n$. Problem (50) can be given by

$$\min \|((\boldsymbol{x} - \overline{\boldsymbol{x}})^T, vec(\boldsymbol{y} - \overline{\boldsymbol{y}})^T, (\boldsymbol{\mu} - \overline{\boldsymbol{\mu}})^T, (\boldsymbol{\rho} - \overline{\boldsymbol{\rho}})^T, (\boldsymbol{\sigma} - \overline{\boldsymbol{\sigma}})^T, vec(\boldsymbol{u} - \overline{\boldsymbol{u}})^T, (\boldsymbol{t} - \overline{\boldsymbol{t}})^T)^T\|_2$$

s.t.

$$x_i = x_i'^{k0} + x_i'^{k1} \qquad i \in I$$

$$y_{ij} = y_{ij}'^{k0} + y_{ij}'^{k1} \qquad i \in I, j \in J$$

$$\mu_i = \mu_i'^{k0} + \mu_i'^{k1} \qquad i \in I$$

$$\rho_i = \rho_i'^{k0} + \rho_i'^{k1} \qquad i \in I$$

$$\sigma_i = \sigma_i'^{k0} + \sigma_i'^{k1} \qquad i \in I$$

$$u_{ij} = u_{ij}'^{k0} + u_{ij}'^{k1} \qquad i \in I, j \in J$$

$$t_i = t_i'^{k0} + t_i'^{k1} \qquad i \in I$$

$$\sum_{j \in J} \lambda_j y_{ij}'^{k0} < \mu_i'^{k0} \qquad i \in I$$

$$\sum_{j \in J} \lambda_j y_{ij}'^{k1} < \mu_i'^{k1} \qquad i \in I$$

$$m_i x_i'^{k0} \leq \mu_i'^{k0} \leq M_i x_i'^{k0} \qquad i \in I$$

$$m_i x_i'^{k1} \leq \mu_i'^{k1} \leq M_i x_i'^{k1} \qquad i \in I$$

$$\sigma_i'^{k0} - (\lambda_{k0} - y_{ij}'^{k0}) Q_i \leq u_{ij}'^{k0} \leq \sigma_i'^{k0} \qquad i \in I, j \in J$$

$$\sigma_i'^{k1} - (\lambda_{k1} - y_{ij}'^{k1}) Q_i \leq u_{ij}'^{k1} \leq \sigma_i'^{k1} \qquad i \in I, j \in J$$

$$0 \leq u_{ij}'^{k0} \leq y_{ij}'^{k0} Q_i \qquad i \in I, j \in J$$

$$0 \leq u_{ij}'^{k1} \leq y_{ij}'^{k1} Q_i \qquad i \in I, j \in J$$

$$\left\| \begin{array}{c} 2\sqrt{\lambda_1} y_{i1}'^{k0} \\ \vdots \\ 2\sqrt{\lambda_{|J|}} y_{i|J|}'^{k0} \\ \rho_i'^{k0} - \mu_i'^{k0} \end{array} \right\| \leq \rho_i'^{k0} + \mu_i'^{k0} \qquad i \in I$$

$$\left\| \begin{array}{c} 2\sqrt{\lambda_1} y_{i1}'^{k1} \\ \vdots \\ 2\sqrt{\lambda_{|J|}} y_{i|J|}'^{k1} \\ \rho_i'^{k1} - \mu_i'^{k1} \end{array} \right\| \leq \rho_i'^{k1} + \mu_i'^{k1} \qquad i \in I$$

$$\left\| \begin{array}{c} 2\rho_i'^{k0} \\ 2\sum_{j\in J} \lambda_j u_{ij}'^{k0} \\ 2\lambda_{k0} - 2\rho_i'^{k0} - t_i'^{k0} \end{array} \right\| \leq 2\lambda_{k0} - 2\rho_i'^{k0} + t_i'^{k0} \qquad i \in I$$

$$\left\| \begin{array}{c} 2\rho_i'^{k1} \\ 2\sum_{j\in J} \lambda_j u_{ij}'^{k1} \\ 2\lambda_{k1} - 2\rho_i'^{k1} - t_i'^{k1} \end{array} \right\| \leq 2\lambda_{k1} - 2\rho_i'^{k1} + t_i'^{k1} \qquad i \in I$$

$$t_i'^{k0} \geq 0, \rho_i'^{k0} \geq 0, \qquad i \in I$$

$$t_i'^{k1} \geq 0, \rho_i'^{k1} \geq 0, \qquad i \in I$$

$$0 \leq x_i'^{k0} \leq \lambda_{k0}, 0 \leq y_{ij}'^{k0} \leq \lambda_{k0} \qquad i \in I, j \in J$$

$$0 \leq x_i'^{k1} \leq \lambda_{k1}, 0 \leq y_{ij}'^{k1} \leq \lambda_{k1} \qquad i \in I, j \in J$$

$$\left(x'^{k0}\right)_k = 0, \left(x'^{k1}\right)_k = \lambda_{k1}$$

$$\lambda_{k0} + \lambda_{k1} = 1, \quad \lambda_{k0}, \lambda_{k1} \geq 0$$

$$v_i\left(\frac{\mu_i^{k0}}{\lambda_{k0}}\right) \leq \left(\frac{\sigma_i^{k0}}{\lambda_{k0}}\right)^2 \qquad i \in I$$

$$v_i\left(\frac{\mu_i^{k1}}{\lambda_{k1}}\right) \leq \left(\frac{\sigma_i^{k1}}{\lambda_{k1}}\right)^2 \qquad i \in I.$$

According to Proposition 1, the last two sets of constraints are SOC-representable as long as the service times have the stochastic representation (24). Therefore, considering the quadratic objective function, the minimum distance problem is transformable to an SOCP.

### Appendix C. An alternative MISOCP reformulation

Our problem may be formulated in other ways as MISOCP. This appendix provides one of these alternatives; its adaptations for the special cases in Section 3 can similarly be derived. Our numerical study shows that this formulation has far less computational efficiency and it is why we have not used it in the paper. However, the proposed valid inequalities are designed based on the same idea used in this reformulation.

By substituting $\rho_i = \frac{\sum_{j \in J} \lambda_j y_{ij}}{\mu_i}$, the objective (1) can be formulated as:

$$\min \sum_{i \in I} \left( ec_i x_i + sc_i \mu_i + wc_i \left( \frac{\rho_i^2}{2(1-\rho_i)} + \frac{\left(\sum_{j \in J} \lambda_j y_{ij}\right)^2 \sigma_i^2}{2(1-\rho_i)} + \rho_i \right) + \sum_{j \in J} tc_{ij} \lambda_j y_{ij} \right).$$

By defining new variables $r_i, s_i \geq 0$ for $i \in I$ such that

$$\rho_i^2 = 2(1-\rho_i)t_i \qquad i \in I$$

$$\left(\sum_{j\in J}\lambda_j u_{ij}\right)^2 = 2(1-\rho_i)t_i \qquad i \in I$$

where

$$u_{ij} = \sigma_i y_{ij} \qquad i \in I, j \in J,$$

the model can be rewritten as follows:

$$\min \sum_{i\in I}\left(ec_i x_i + sc_i \mu_i + wc_i(r_i + s_i + \rho_i) + \sum_{j\in J} tc_{ij}\lambda_j y_{ij}\right)$$

s.t.

$$\sum_{i\in I} y_{ij} = 1 \qquad j \in J$$

$$y_{ij} \leq x_i \qquad i \in I, j \in J$$

$$\sum_{j\in J}\lambda_j y_{ij} \leq \mu_i \qquad i \in I$$

$$m_i x_i \leq \mu_i \leq M_i x_i \qquad i \in I$$

$$v_i(\mu_i) \leq \sigma_i^2 \qquad i \in I$$

$$\sigma_i - (1 - y_{ij})Q_i \leq u_{ij} \leq \sigma_i \qquad i \in I, j \in J$$

$$0 \leq u_{ij} \leq y_{ij} Q_i \qquad i \in I, j \in J$$

$$\left\| \begin{array}{c} 2\sqrt{\lambda_1} y_{i1} \\ \vdots \\ 2\sqrt{\lambda_{|J|}} y_{i|J|} \\ \rho_i - \mu_i \end{array} \right\| \leq \rho_i + \mu_i \qquad i \in I$$

$$\left\| \begin{array}{c} 2\rho_i \\ 2 - 2\rho_i - r_i \end{array} \right\| \leq 2 - 2\rho_i + r_i \qquad i \in I$$

$$\left\| \begin{array}{c} 2\sum_{j\in J}\lambda_j u_{ij} \\ 2 - 2\rho_i - s_i \end{array} \right\| \leq 2 - 2\rho_i + s_i \qquad i \in I$$

$$x_i \in \{0,1\}, y_{ij} \in \{0,1\} \qquad i \in I, j \in J$$

$$r_i \geq 0, s_i \geq 0, \rho_i \geq 0 \qquad i \in I.$$